\documentclass[11pt]{article}
\usepackage{epsfig}

\def\be{\begin{equation}}
\def\ee{\end{equation}}

\def\nn{\nonumber}
\def\<{\langle}
\def\>{\rangle}
\def\lb{\label}

\def\R{{\bf R}}

\def\Z{{\bf Z}}
\def\N{{\bf N}}

\def\diag{{\rm diag}}

\def\hb{\vrule height0.18cm width0.14cm $\,$}

\title{Index theory for linear self-adjoint operator equations and
nontrivial \\solutions for asymptotically linear operator
equations(II)\thanks{Partially supported by the National Natural
Science Foundation of China(10871095)}
 }
\author{Yujun Dong\\
Department of Mathematics, Nanjing Normal University, Nanjing,\\
Jiangsu 210097,
P. R. China. E-mail: yjdong@njnu.edu.cn\\
Yuan Shan\\
Department of Mathematics, Nanjing Normal University, Nanjing,\\
Jiangsu 210097, P. R. China. E-mail:yuanfang8887js@163.com}
\date{}
\begin{document}

\maketitle

\begin{abstract}
{\it Reference [1] established an index theory for a class of linear
selfadjoint operator equations covering both second order linear
Hamiltonian systems and first order linear Hamiltonian systems as
special cases. In this paper based upon this index theory we
construct a new reduced functional to investigate multiple solutions
for asymptotically linear operator equations by Morse theory. The
functional is defined on an infinite dimensional Hilbert space, is
twice differentiable and has a finite Morse index. Investigating
critical points of this functional by Morse theory gives us a
unified way to deal with nontrivial solutions of both asymptotically
second order Hamiltonian systems and asymptotically first order
Hamiltonian systems.}
\end{abstract}

\bigskip


\renewcommand{\theequation}{\thesection.\arabic{equation}}

\setcounter{section}{0}
\setcounter{equation}{0}

\section{Introduction and main results}

Let $X$ be a real separable infinite dimensional Hilbert space with
inner product $(\cdot, \cdot )$ and norm $\| \cdot\|$. Let $A:
D(A)\subset X\to X$ be a unbounded linear self-adjoint operator with
domain $D(A)$ satisfying $\sigma(A)=\sigma_d(A)$. We investigate the
following equation
   \be
   Ax-\Phi'(x)=\theta, \lb{1.1}
   \ee
where $\Phi\in C^1(X)$ and $\Phi'(x)$ is the derivative of $\Phi$
with respect to $x$ in $X$. This equation covers both second order
Hamiltonian systems and first order Hamiltonian systems as special
cases. We will construct a new reduced functional to investigate
(1.1) by Morse theory. The functional is defined on an infinite
subspace of $X $, is twice differentiable and has a finite
dimensional Morse index at its any critical point. The main result
is the following theorem.

{\bf Theorem 1.1}. Assume that

(i) $\Phi''(x)$ exists and is bounded for $x\in X$,
$\Phi'(\theta)=\theta$, $\Phi\in C^{2}(V)$ with
$V:=\hbox{D}(|A|^{\frac{1}{2}})$;

(ii) there exists $B_{1}$, $B_{2}\in L_{s}(X)$ satisfying
$i_{A}(B_{1})=i_{A}(B_{2})$, $\nu_{A}(B_{2})=0$ and $B:X \to L_s(X),
C: X\to X$ such that
   \begin{eqnarray}
&&\Phi'(x)=B(x)x+C(x)\ \  \hbox{for any} \ \  x,\nn\\
&& B_1\leq B(x)\leq B_2,\ \ C(x)\ \ \hbox{is bounded};\nn
   \end{eqnarray}

(iii)with $B_{0}:=\Phi''(\theta)$ we have
  $$
i_{A}(B_{1})\notin [i_{A}(B_{0}),i_{A}(B_{0})+\nu_{A}(B_{0}) ].
  $$
Then (1.1) has a nontrivial solution $x=x_0$.

Under the further assumption that

(iv) $\nu_{A}(B_{0})=0$ and $|i_{A}(B_{1})-i_{A}(B_{0})|\geq
\nu_{A}(\Phi''(x_0))$, (1.1) has two nontrivial solutions.

\vskip4mm

In the theorem we used notations $(i_A(B),\nu_A(B))$ concerning the
linear selfadjoint operator equation
 \be
 Ax-Bx=0\lb{1.2}
 \ee
for any $B\in L_s(X)$, which will be defined as follows.

{\bf Definition 1.2} For any $B\in  L_s(X)$, we define
  \be
  \nu_A(B)=\dim \ker(A-B).
  \ee
$\nu_A(B)$ is called the nullity of $B$.

{\bf Definition 1.3} For any $B_1,B_2\in L_s(X)$ with $B_1<B_2$, we
define
  \be
  I_A(B_1,B_2)=\sum_{\lambda\in [0,1)}\nu_A((1-\lambda)B_1+\lambda
  B_2);\lb{1.3}
  \ee
and for any $B_1,B_2\in L_s(X)$ we define
  \be
  I_A(B_1,B_2)=I_A(B_1,k I)-I_A(B_2,k I)\lb{1.4}
  \ee
where $I:X\to X$ is the identity map and $k I>B_1, k I>B_2$ for some
real number $k>0$. We call $I_A(B_1,B_2)$ the relative Morse index
between $B_1$ and $B_2$.

\vskip 2mm Let $B_0\in L_s(X)$ be fixed and let $i_A(B_0)$ be a
prescribed integer associated with $B_0$.

{\bf Definition 1.4} For any  $B\in L_s(X)$ we define
   \be
   i_A(B)=i_A(B_0)+I_A(B_0,B).\lb{1.5}
   \ee

\vskip2mm As in [1] we call $i_A(B)$ the index of $B$ and $i_A(B_0)$
is called initial index. Generally, the initial index can be any
prescribed integer and the index $i_A(B)$ also depends on $B_0$ and
the initial index. Let $X_1$ be a nontrivial subspace of $X$. For
$B_1,B_2\in {\cal L}_s(X)$ we write $B_1\leq B_2$ with respect to
$X_1$ if and only if $(B_1x,x)\leq (B_2x,x)$ for any $x\in X_1$; we
write $B_1<B_2$ with respect to $X_1$ if and only if
$(B_1x,x)<(B_2x,x)$ for any $x\in X_1\backslash\{\theta\}$. If
$X_1=X$ we just write $B_1\leq B_2$ or $B_1<B_2$.

\vskip2mm {\bf Theorem 1.5}  (i) For any $B,B_1, B_2\in L_s(X)$,
$\nu_A(B)\in \N, I_A(B_1,B_2)\in \Z$ and $i_A(B)\in \Z$ are
well-defined;

  (ii) For any $B_1,B_2\in L_s(X)$,
  $I_A(B_1,B_2)=i_A(B_2)-i_A(B_1)$, and if $B_{1}<B_{2}$ with
  respect to Ker$(A-(1-\lambda)B_{1}-\lambda B_{2})\neq \{
  \theta\}$ for $t\in [0,1)$, then (1.4) holds;

  (iii) For any $B_1,B_2, B_3\in L_s(X)$,
   $I_A(B_1,B_2)+I_A(B_2,B_3)=I_A(B_1,B_3)$;

  (iv) For any $B_1,B_2\in L_s(X)$, if $B_1\leq B_2$,
then $i_A(B_1)\leq i_A(B_2), \nu_A(B_1)+i_A(B_1)\leq\nu_A(B_2)+
i_A(B_2)$; if $B_1<B_2$ with respect to Ker$(A-B_{1})$, then
$\nu_A(B_1)+i_A(B_1)\leq i_A(B_2)$.

  (v) If there exists $B_0\in  L_s(X)$ such
that $\sum_{\lambda<0}\nu_A(B_0+\lambda I)<+\infty$, we will choose
this integer for $i_A(B_0)$. Then the index defined by Definition
1.3 satisfies
  \be
  i_A(B)=\sum_{\lambda<0}\nu_A(B+\lambda I).\lb{1.6}
  \ee

\vskip2mm In [1] an index theory for $Ax+Bx=\theta$ was established
by the concept of relative Morse index and dual variational methods.
Here we discuss $(1.2)$ instead only because the new form will bring
convenience to the proof of Theorem 1.1 as will be seen in Sections
4-5. In [1] it was assumed that $A$ satisfies the following
condition:

   {\it $(A)$: $A: Y\to X$ is linear bounded, symmetric i.e.
    $(Ax,y)=(x,Ay)$ for any $x,y\in Y$, $R(A)$ is closed in $X$ and
    $X=R(A)\bigoplus ker(A)$, where $X$ is a real separable infinite
    dimensional Hilbert space, $Y\subset X$ is a Banach space and the embedding
     $Y\hookrightarrow X$ is compact.}

In order to prove Theorem 1.5 we first prove the following
proposition.

{\bf Proposition 1.6} $A: D(A)\subset X\to X$ is selfadjoint and
$\sigma(A)=\sigma_d(A)$ if and only if $A$ satisfies condition
$(A)$.

 {\bf Proof}. Sufficiency:  Denote all the eigenvalues of $A$ with multiplicities by
$\{\lambda_j\}_{j=-\infty}^{\infty}$ satisfying
$\lambda_j\leq\lambda_{j+1} \forall j$ and $\lambda_j\to\pm\infty$
as $j\to\pm\infty$. There is a unit orthogonal basis $\{e_j\}$ of
$X$ such that $X=\{\sum_{j=-\infty}^{\infty} c_je_j|\sum
c_j^2<\infty\}$ and $D(A)=\{\sum_{j=-\infty}^{\infty}
c_je_j|\sum_{j=-\infty}^{\infty}(1+\lambda_j^2)c_j^2<\infty\}$. For
any $x=\sum_{j=-\infty}^{\infty} c_je_j\in X$, because
$x_n:=\sum_{j=-n}^n c_je_j\in D(A)$ and $x_n\to x$ in $X$. Thus,
$D(A)$ is dense in $X$. In order to prove the adjointness of $A$ we
need only $D(A^{\ast})\subset D(A)$. In fact, assume
$x=\sum_{j=-\infty}^{\infty}c_je_j\in D(A^{\ast})$. By definition
there exists a constant $C>0$ such that $|(Ay,x)|\leq C||y||$ for
any $y\in D(A)$. If we choose $y=\sum_{j=-n}^n\lambda_jc_je_j\in
D(A)\ \forall \ n$, then $\sum_{j=-n}^n\lambda_j^2c_j^2\leq C^2$.
Thus $x\in D(A)$.

Necessity: Because $\sigma(A)=\sigma_d(A)$, by definition,
$D(A)=\{x\in X|\Sigma_{\lambda\in
\sigma_d(A)}(1+\lambda^2)\|E(\{\lambda\})x\|^2<\infty\}$; and
$\forall  \lambda\in \sigma_d(A)$, $\lambda$ is isolated and
$E(\{\lambda\})X=ker(A-\lambda I)$, the embedding $(D(A),
\|\cdot\|_G)\hookrightarrow X$ is compact. To finish the proof, we
prove $R(A)$ is closed. In fact, any $x\in X$ satisfies
$x=\sum_{\lambda\in \sigma_d(A)}E(\{\lambda\})x$. If $x\perp
ker(A)$, then $x=\sum_{\lambda\in
\sigma_d(A)\backslash\{0\}}E(\{\lambda\})x$ and $\sum_{\lambda\in
\sigma_d(A)\backslash\{0\}}{1\over \lambda}E(\{\lambda\})x\in D(A)$.
This means that $R(A)=(ker A)^{\perp}$ is closed.\hfill\hb

{\bf Proof of Theorem 1.5} We give two proofs.

Step 1: Set $A_1:=-A$. Because $A$ is selfadjoint and
$\sigma(A)=\sigma_d(A)$, so does $A_1$. By Proposition 1.6, $A_1$
satisfies condition $(A)$. From [1, Definitions 3.1.1, 3.1.2 and
3.1.3] $(i_{A_1}(B),\nu_{A_1}(B))$ is defined. Denote
$(i_{A_1}(B),\nu_{A_1}(B))$ by $(i_A(B),\nu_A(B))$. Then (1.3),
(1.4), (1.5) and (1.6) are satisfied. And [1, Propositions 3.1.4 and
3.1.5, and Lemma 3.2.1] imply all the conclusions of Theorem 1.4.

Step 2: Because $\sigma(A)=\sigma_{d}(A)$, every eigenvalue is
isolated and corresponds to a finite dimensional subspace of
eigenvectors. Let $-\mu\notin \sigma(A)$. Then $A+\mu I$ is
invertible. For any $B\in L_s(X)$ and such a $\mu>0$ large enough
satisfying $B+\mu I>I$, the following bilinear form
   $$
\psi_{A,\mu;B}(x,y)=((B+\mu I)^{-1}x,y)-((A+\mu I)^{-1}x,y), \ \
\forall x,y\in X
  $$
has finite Morse index $m^-_{A,\mu}(B)$ and finite Morse nullity
$m^0_{A,\mu}(B)$. Because the first term in $\psi_{A,\mu;B}$ is the
same with the second term in $\psi_{A,B|B_0}$ defined by (3.5) in
[1] when $B_0=\mu I$, and the left terms in both these two bilinear
forms do not depend on $B$, all the associated conclusions of [1,
Theorem 3.2.4] hold, and from which in a way similar to the proofs
of [1, Propositions 3.1.4 and 3.1.5] we can complete the proof.
\hfill\hb

\vskip2mm Many authors investigated (1.1) and its special
cases-first order or second order Hamiltonian systems. We refer to
[2,3,4,5,6,7,8,17] for references. Reference [1] investigated (1.1)
by Morse theory only in the case $\sigma(A)$ is bounded from below
and the result cannot apply to first order Hamiltonian systems.
However, Theorem 1.1 applies to both second order Hamiltonian
systems and first order Hamiltonian systems respectively. The paper
will be organized in the following way. In Section 2 as applications
of Theorem 1.1 we investigate first order and second order
asymptotically Hamiltonian systems. In Section 3, we construct the
mentioned reduced functional. In Section 4 we investigate properties
of the Morse index of the functional at any critical point. And in
the last section we prove Theorem 1.1 by Morse theory.

\setcounter{equation}{0}
\section{Applications of Theorem 1.1:
          Nontrivial solutions for asymptotically linear Hamiltonian systems}

In this section as applications of Theorem 1.1 we investigate
nontrivial solutions of both first order Hamiltonian systems and
second order Hamiltonian systems satisfying various boundary value
conditions.

 \vskip4mm\noindent{\bf 2.1 $\quad$ Hamiltonian systems
satisfying Bolza boundary value conditions}

As in [1, Section 3.4] we are interested in the index theory for the
following Hamiltonian system
    \begin{eqnarray}
  &&-J\dot{x}-B(t)x=0\lb{2.1}\\
  &&x_1(0)\cos\alpha+x_2(0)\sin\alpha=0\lb{2.2}\\
  &&x_1(1)\cos\beta+x_2(1)\sin\beta=0\lb{2.3}
 \end{eqnarray}
where $B\in L^{\infty}((0,1);GL_s(\R^{2n})), 0\leq\alpha<\pi$ and
$0<\beta\leq\pi, x=(x_1,x_2)\in\R^n\times\R^n$. Define
$X:=L^2([0,1];\R^{2n})$, $D(A)=\{x\in H^1([0,1];\R^{2n})|x$
satisfies $(2.2-2.3)\}$, $(Ax)(t)=-J\dot{x}(t)$ for any $x\in D(A)$,
and $(\bar{B}x)(t)=B(x)x(t)$ for $x\in L^{2}([0,1],\textbf{R}^{2n})$
. Then $D(A)$ with the graph norm $\|\cdot\|_G$ is a Banach space
and the embedding from $D(A)$ to $X$ is compact. It is also easy to
check that $A$ is symmetric i.e. $(Ax,y)=(x,Ay)$ for any $x,y\in
D(A)$. As proved in [1] there hold that $R(A)$ is closed in $X$ and
$X=R(A)\oplus $ ker$(A)$. So we have the following definitions and
properties.

\vskip4mm
 {\bf Definition 2.1[1, Definition 3.4.4]} For any $B\in L^{\infty}((0,1);GL_s(\R^{2n}))$,
we define
  \begin{eqnarray}
&&\nu^f_{\alpha,\beta}(B):=\dim \ker(A-\bar{B}),\nn\\
&&i^f_{\alpha,\beta}(\diag\{0,I_n\}):=i^s_{I_n,\alpha,\beta}(0),\nn\\
&&i^f_{\alpha,\beta}(B):=i^f_{\alpha,\beta}(\diag\{0,I_n\})+I^f_{\alpha,\beta}(\diag\{0,I_n\},B);\nn
  \end{eqnarray}
and
  \begin{eqnarray}
  &&I^f_{\alpha,\beta}(B_1,B_2)=\sum_{\lambda\in [0,1)}\nu^f_{\alpha,\beta}((1-\lambda)B_1+\lambda B_2)
  \,\,\,\,\hbox{as}\,\,\,\,B_1<B_2,\nn\\
  &&I^f_{\alpha,\beta}(B_1,B_2)=I^f_{\alpha,\beta}(B_1,k I)-I^f_{\alpha,\beta}(B_2,kI)
          \,\,\,\,\hbox{for every}\,\,\,\, B_1,B_2\,\,\,\,\hbox{with}\,\,\,\,k I >B_1, k I>B_2\nn
  \end{eqnarray}
where for any  $B\in L^{\infty}([0,1];GL_s(\R^n))$,
$(i^s_{I_n,\alpha,\beta}(B),\nu^s_{I_n,\alpha,\beta}(B))\in\N\times\{0,1,2,\cdots,n\}$
will be defined in Definition 2.5.

\vskip4mm

{\bf Proposition 2.2[1, Proposition 3.4.2]} We have the following
properties:

 (i) For any $B\in L^{\infty}((0,1);GL_s(\R^{2n}))$, $\nu^f_{\alpha,\beta}(B)$ is
the dimension of the solution subspace of system (2.1-2.3) and
  $$
 (i^f_{\alpha,\beta}(B),\nu^f_{I_n,\alpha,\beta}(B))\in\Z\times\{0,1,2,\cdots,n\}.
  $$

 (ii) For any $B_1, B_2\in L^{\infty}((0,1);$GL$_s(\R^{2n}))$, if
$B_1\leq B_2$, then $i^f_{\alpha,\beta}(B_1)\leq
i^f_{\alpha,\beta}(B_2)$ and
$i^f_{\alpha,\beta}(B_1)+\nu^f_{\alpha,\beta}(B_1)\leq
i^f_{\alpha,\beta}(B_2)+\nu^f_{\alpha,\beta}(B_1)$; if $B_1<B_2$
then $i^f_{\alpha,\beta}(B_1)+\nu^f_{\alpha,\beta}(B_1)\leq
i^f_{\alpha,\beta}(B_2)$.

 (iii) For any $B\in L^{\infty}((0,1);$GL$_s(\R^n))$, there holds
   $$
   (i^f_{\alpha,\beta}(\diag\{B,I_n\}),\nu^f_{\alpha,\beta}(\diag\{B,I_n\}))
   =(i^s_{I_n,\alpha,\beta}(B),\nu^s_{I_n,\alpha,\beta}(B)).
   $$

   Here for any $B_{1}$, $B_{2}\in
   L^{\infty}([0,1];GL_s(\R^{2n}))$, we write $B_{1}\leq B_{2}$ if
   $B_{1}(t)\leq B_{2}(t)$ for a.e. $t\in[0,1]$; and we write
   $B_{1}<B_{2}$ if $B_{1}\leq B_{2}$ and $B_{1}(t)< B_{2}(t)$ for
   $t$ belonging to a subset of $(0,1)$ with nonzero measure. If $B_{1}\leq
   B_{2}$ then $\bar{B}_{1}\leq \bar{B}_{2}$; and if $B_{1}<B_{2}$
   then $\bar{B}_{1}<\bar{B}_{2}$ with respect to
   Ker$(A-\bar{B}_{1})$. So Proposition 2.2(ii) follows from
   Theorem1.5(iv) directly.

We now use this index theory to discuss the solvability of the
following Hamiltonian system (2.2-2.3) and
 \be
  -J\dot{x}- H'(t,x)=\theta,\lb{2.4}
    \ee
where $H:[0,1]\times\R^{2n}\to \R^{2n}$ is differentiable and
$H'(t,x)$ is the gradient of $H$ with respect to $x$.

Let $H''(t,x)$ denote the second derivative of $H(t,x)$ with respect
to $x$. We have the following theorem.

 {\bf Theorem 2.3}. Assume that

(i) $H''(t,x)$ is continuous and is bounded for $(t,x)\in
[0,1]\times \R^{2n}$, and $H'(t,\theta)=\theta$;

(ii) there exists $B_{1}$, $B_{2}\in
L^{\infty}([0,1];GL_{s}(\R^{2n}))$ satisfying
$i^f_{\alpha,\beta}(B_{1})=i^f_{\alpha,\beta}(B_{2})$,
$\nu^f_{\alpha,\beta}(B_{2})=0$ and
  $$
B_1(t)\leq H''(t,x)\leq B_2(t)\ \ \forall (t,x)\in [0,1]\times
\R^{2n}\ \ with \ \ |x|\geq r>0;
 $$

(iii)with $B_{0}:=H''(\cdot,\theta)$ we have
   $$
i^f_{\alpha,\beta}(B_{1})\notin
[i^f_{\alpha,\beta}(B_{0}),i^f_{\alpha,\beta}(B_{0})+\nu^f_{\alpha,\beta}(B_{0})
].
    $$
Then (2.4)(2.2-2.3) has a nontrivial solution.

Under the further assumption that

(iv) $\nu^f_{\alpha,\beta}(B_{0})=0$ and
$|i^f_{\alpha,\beta}(B_{1})-i^f_{\alpha,\beta}(B_{0})|\geq n $,
(2.4)(2.2-2.3) has two nontrivial solutions.

 Define
  \begin{eqnarray}
 \Phi(x)=\int_0^1H(t,x(t))dt \ \ \forall  x\in X
 \end{eqnarray}

It is easy to check that (2.4)(2.2-2.3) is equivalent to (1.1). So
in order to prove Theorem 2.3 by Theorem 1.1 we only need the
following lemma.

{\bf Lemma 2.4} (i) $\Phi\in C^2(V)$;

(ii) Under assumption(ii) of Theorem 2.3, there holds
 \begin{eqnarray}
H'(t,x)=B(t,x)x+C(t,x)\nn\\
B_{1}(t)-\epsilon_1 I_{2n}\leq B(t,x)\leq B_{2}(t)+\epsilon_1
I_{2n}\nn
 \end{eqnarray}
where for any $x\in X$, $B(\cdot,x(\cdot))\in L_s(X)$ $\epsilon_1>0$
satisfying $i^f_{\alpha,\beta}(B_{1}-\epsilon_1
I_{2n})=i^f_{\alpha,\beta}(B_{1})=i^f_{\alpha,\beta}(B_{2})$,
$\nu^f_{\alpha,\beta}(B_{2}+\epsilon_1 I_{2n})=0$, and
$C(\cdot,x(\cdot))\in X$ is uniformly bounded.

{\bf Proof}. (i) From [1] it follows that
$\sigma(A)=\sigma_d(A)\subset \R$, and $\sigma_d(A)$ is unbounded
from both above and from bellow. In fact, if $\lambda\neq 0$ is an
eigenvalue of $A$ with an eigenvector $e^{\lambda J t}c$ for
$c\in\R^{2n}$, then for any integer $k$, $\lambda+2k\pi$ is also an
eigenvalue with the eigenvector $e^{(\lambda +2\pi)J t}c$. Note that
$V=\{\sum_{j=-\infty}^{\infty}
c_je_j|\sum_{j=-\infty}^{\infty}(1+|\lambda_j|)c_j^2<\infty\}$,
where $\{e_{j} \}$ is an orthonormal basis of $X$, $\{\lambda_{j}
\}$ is all the eigenvalues of $A$ with multiplicities and
$\lambda_{j}\leq \lambda_{j+1}$ $\forall j$. For any $x,x_{0} \in V$
and $u=\sum_{j=-\infty}^{\infty} c_je_j, v=\sum_{j=-\infty}^{\infty}
c^{'}_je_j\in V$ satisfying $||u||_V\leq 1$ and $||v||_V\leq 1$, we
have $|c_j |\leq 1$, $|c^{'}_j |\leq 1$, and by assumption (ii)
$\exists M>0$ such that $||H''(\cdot,x(\cdot))u||\leq M||u||$ for
any $x,u\in X$. Hence,
    \begin{eqnarray}
&&||\Phi''_V(x)-\Phi''_V(x_0)||_V=\sup_{||u||_V\leq 1,||v||_V\leq 1}|
                                         \int_0^1((H''(t,x(t))-H''(t,x_0(t)))u(t),v(t))dt|\nn\\
&&\ \
\leq\sum_{i,j=-n}^n|\int_0^1((H''(t,x(t))-H''(t,x_0(t)))e_j(t),e_i(t))dt|
    +4M({1\over |\lambda_{-n-1}|}+{1\over |\lambda_{n+1}|})\nn
   \end{eqnarray}
Because $\lambda_j\to\pm\infty$ as $j\to\pm\infty$ and from Theorem
4 in page 97 of Ekeland's book[13] it follows that for fixed $i,j$
as $x\to x_0$ in $X$
  $$
\int_0^1((H''(t,x(t))-H''(t,x_0(t)))e_j(t),e_i(t))dt\to 0;
  $$
we obtain
   $$
||\Phi''_V(x)-\Phi''_V(x_0)||_V\to 0
  $$
as $x\to x_0$ in $V$.

(ii) From assumption (ii) there exists $\epsilon_1>0$ such that
$i_{\alpha,\beta}^f(B_{1}-\epsilon_1
I_n)=i_{\alpha,\beta}^f(B_{1})$,
$i_{\alpha,\beta}^f(B_{2}+\epsilon_1I_n)+\nu_{\alpha,\beta}^f(B_{2}+\epsilon_1I_n)=i_{\alpha,\beta}^f(B_{2})$.
And we can choose $\delta \in (0,1)$ such that
   \begin{eqnarray}
-\epsilon_1/2 \leq \delta B_{1}\leq \delta B_{2} \leq \epsilon_1/2,
\ \ \hbox{ and }\ \  \delta M< \frac{1}{2}\epsilon_1.\nn
   \end{eqnarray}
Define
   \begin{eqnarray}
B(t,x)&&= \int_0^1H''(t,\theta x)d\theta\ \  \hbox{ if } |x|\geq r/\delta,\nn\\
    &&=B_{1}(t)   \hbox{    otherwise};\nn
  \end{eqnarray}
and $C(t,x)=H'(t,x)-B(t,x)x$. Then for any $x\in X, B(\cdot,
x(\cdot))\in L_s(X), B_1-\epsilon I_{2n}\leq B(\cdot,
x(\cdot))\leq B_2+\epsilon I_{2n}$ and $C(\cdot,x(\cdot))\in X$ is
bounded uniformly for $x\in X$. The proof is complete.  \hfill\hb

\vskip2mm {\bf Remark}. In the special case $\alpha=0,\beta=\pi$,
Theorem 2.3 was given in [10]. However the proof there is not
correct because generally the integral functional defined by (2.5)
is not twice differentiable in $L^2([0,1];\R^{2n})$ even we assume
$H''(t,x)$ is continuous and bounded for $(t,x)\in
[0,1]\times\R^{2n}$.

{\vskip2mm} Theorem 1.1 can also be used to investigate second order
Hamiltonian systems satisfying Sturm-Liouville boundary value
conditions. Recall that an index theory has been established in [1]
for the following system:
    \begin{eqnarray}
  &&-\ddot{x}-B(t)x=0\lb{2.5}\\
  &&x(0)\cos\alpha-x'(0)\sin\alpha=0\lb{2.6}\\
  &&x(1)\cos\beta-x'(1)\sin\beta=0\lb{2.7}
    \end{eqnarray}
where $0\leq\alpha<\pi$ and $0<\beta\leq\pi$.

{\bf Definition 2.5[1, Definition 2.3.2 and Proposition 2.3.3]} We
define
   \begin{eqnarray}
&&\nu^s_{\alpha,\beta}(B)\ \ \hbox{is the dimension of the solution subspace of} \ \ (2.6-2.8)\ \   \nn\\
&&i^s_{\alpha,\beta}(B):=\sum_{\lambda<0}\nu^s_{\alpha,\beta}(B+\lambda
I_n).\nn
   \end{eqnarray}
As before $(i^s_{\alpha,\beta}(B),\nu^s_{\alpha,\beta}(B))$ has
useful properties, which can be found in [1]. This index can be
used to investigate the following nonlinear Hamiltonian system
(2.7-2.8) and
 \be
  -\ddot{x}-V'(t,x)=0,\lb{2.8}
  \ee
where $V:[0,1]\times\R^n\to \R^n$ is continuous and $V'(t,x)$
denotes the gradient of $V(t,x)$ with respect to $x$. Let $V''(t,x)$
denote the second derivative of $V(t,x)$ with respect to $x$. We
have the following theorem.

{\bf Theorem 2.6} Assume that

(i) $V''(t,x)$ is continuous and is bounded for $(t,x)\in
[0,1]\times \R^{2n}$, and $H'(t,\theta)=\theta$;

(ii) there exist $B_{1}$, $B_{2}\in
L^{\infty}([0,1];GL_{s}(\R^{n}))$ satisfying
$i^s_{\alpha,\beta}(B_{1})=i^s_{\alpha,\beta}(B_{2})$,
$\nu^s_{\alpha,\beta}(B_{2})=0$ and
  $$
B_1(t)\leq V''(t,x)\leq B_2(t)\ \ \forall (t,x)\in [0,1]\times
\R^{n}\ \ with \ \ |x|\geq r>0;
 $$

(iii)with $B_{0}:=V''(\cdot,\theta)$ we have
  $$
i^s_{\alpha,\beta}(B_{1})\notin
[i^s_{\alpha,\beta}(B_{0}),i^s_{\alpha,\beta}(B_{0})+\nu^s_{\alpha,\beta}(B_{0})
].
   $$
Then (2.9)(2.7-2.8) has a nontrivial solution.

Under the further assumption that

(iv) $\nu^s_{\alpha,\beta}(B_{0})=0$ and
$|i^s_{\alpha,\beta}(B_{1})-i^s_{\alpha,\beta}(B_{0})|\geq n $,
(2.9)(2.7-2.8) has two nontrivial solutions.

{\bf Proof }  Define $x_1=x, x_2=-\dot{x}, x=(x_1,x_2)$ and
$H(t,x)=V(t,x_1)+{1\over 2}|x_2|^2$. Then (2.9)(2.7-2.8) is
equivalent to (2.4)(2.2-2.3). Under assumption (ii) we have
$V'(t,x_1)=B(t,x_1)x_1+C(t,x_1)$ as before. Because
$i^s_{\alpha,\beta}(B)=i^f_{\alpha,\beta}(\diag\{B,I_n\})$, the
result follows. \hfill\hb

{\bf Remark} For the special case $\alpha=0,\beta=\pi$ this theorem
was obtained in [1, Theorem 2.3.7], and [9, Theorem 3.3] by
different methods.

\vskip4mm\noindent{\bf 2.2 $\quad$ Hamiltonian systems satisfying
periodic boundary value conditions}

\vskip4mm Consider the following linear system
     \begin{eqnarray}
  &&-J\dot{x}-B(t)x=0\nn\\
  &&x(1)=Px(0) \lb{2.9}
    \end{eqnarray}
where $P\in Sp(2n)$ is prescribed. Define $X:=L^2([0,1];\R^{2n}),
D(A):=\{x:\in H^1([0,1];\R^{2n})|x$ satisfies $(2.10)\}$. Then the
embedding from $D(A)$ to $X$ is compact. Define
$(Ax)(t):=-J\dot{x}(t)$ for every $x\in D(A)$. Similar to
Proposition 7 in page 22 of Ekeland's book[13], for the given $P\in
Sp(2n)$ there exists $\lambda\in\R$ such that $(e^{J\lambda}-P)c=0$
for some $c\neq 0$. So $\lambda$ is an eigenvalue of $A$ with an
eigenvector $e^{Jt\lambda}c$. We can check $\lambda+2k\pi$ is also
an eigenvalue of $A$ with the eigenvector $e^{Jt(\lambda+2\pi)}c$.
As in Lemma 2.4 $A$ is selfadjoint and $\sigma(A)=\sigma_d(A)$ is
unbounded from both bellow and above.

Choose $i^f_P(0):=i_P(I_{2n})$ defined by Definition 2.2 in [11]. We
have the following definition.

{\bf Definition 2.7[1, Definition 3.5.1]} For any $B\in
L^{\infty}((0,1);$GL$_s(\R^{2n}))$, we define
  \begin{eqnarray}
  &&\nu^f_P(B)=dim \ker(A-B),\nn\\
  &&i^f_P(B)=i^f_P(0)+I^f_P(0,B);\nn
  \end{eqnarray}
and
  \begin{eqnarray}
  &&I^f_P(B_1,B_2)=\sum_{\lambda\in [0,1)}\nu^f_P((1-\lambda)B_1+\lambda B_2)
  \,\,\,\,\hbox{as}\,\,\,\,B_1<B_2,\nn\\
  &&I^f_P(B_1,B_2)=I^f_P(B_1,k id)-I^f_P(B_2,k id)
          \,\,\,\,\hbox{for every}\,\,\,\,B_1,B_2\,\,\,\,\hbox{with}\,\,\,\, k I >B_1, k I>B_2.\nn
  \end{eqnarray}
From Theorem 1.5 we have the following proposition.

{\bf Proposition 2.8[1, Proposition 3.5.2]}. (i)For any $B\in
L^{\infty}((0,1);$GL$_s(\R^{2n}))$, we have
$\nu^f_P(B)\in\{0,1,2,\cdots,2n\}$.

(ii) For any $B_1, B_2\in L^{\infty}((0,1);$GL$_s(\R^{2n}))$
satisfying $B_1< B_2$, we have $i^f_P(B_1)+\nu^f_P(B_1)\leq
i^f_P(B_2)$.

We now discuss the solvability of the following nonlinear system
(2.4) (2.10)
 \begin{eqnarray}
  &&-J\dot{x}-H'(t,x)=0\nn\\
  &&x(1)=Px(0)\nn
  \end{eqnarray}
where $H:[0,1]\times \R^{2n}\to \R$ is differentiable and $P\in
Sp(2n)$ is prescribed.

 {\bf Theorem 2.9}. Assume that

(i) $H''(t,x)$ is continuous and is bounded for $(t,x)\in
[0,1]\times \R^{2n}$, and $H'(t,\theta)=\theta$;

(ii) there exists $B_{1}$, $B_{2}\in
L^{\infty}([0,1];GL_{s}(\R^{2n}))$ satisfying
$i^f_P(B_{1})=i^f_P(B_{2})$, $\nu^f_P(B_{2})=0$ and
  $$
B_1(t)\leq H''(t,x)\leq B_2(t)\ \ \forall (t,x)\in [0,1]\times
\R^{2n}\ \ with \ \ |x|\geq r>0;
 $$

(iii)with $B_{0}:=H''(\cdot,\theta)$ we have
\begin{eqnarray}
i^f_P(B_{1})\notin  [i^f_P(B_{0}),i^f_P(B_{0})+\nu^f_P(B_{0}) ].\nn
\end{eqnarray}
Then (2.4)(2.10) has a nontrivial solution.

Under the further assumption that

(iii) $\nu^f_P(B_{0})=0$ and $|i^f_P(B_{1})-i^f_P(B_{0})|\geq 2n$, (2.4)(2.9) has two nontrivial solutions.

\vskip2mm Similar to Lemma 2.4 $\Phi\in C^2(V)$. Also similar to
Theorem 2.6, Theorem 2.9 follows from Theorem 1.1.

Set $i^f(B)=i_{I_{2n}}^f(B)$ for any $B\in
L^{\infty}([0,1];$GL$_s(\R^{2n})), i^s(B)=i^f(\{B,I_n\})$ for any
$B\in L^{\infty}([0,1];$GL$_s(\R^{n}))$. Concerning periodic
solutions of Hamiltonian systems we have the following theorems.

 {\bf Theorem 2.10}. Assume that

(i) $H''(t,x)$ is continuous and is bounded for $(t,x)\in
[0,1]\times \R^{2n}$, and $H'(t,\theta)=\theta$;

(ii) there exist $B_{1}$, $B_{2}\in
L^{\infty}([0,1];GL_{s}(\R^{2n}))$ satisfying
$i^f(B_{1})=i^f(B_{2})$, $\nu^f(B_{2})=0$ and
  $$
B_1(t)\leq H''(t,x)\leq B_2(t)\ \ \forall (t,x)\in [0,1]\times
\R^{2n}\ \ with \ \ |x|\geq r>0;
 $$

(iii)with $B_{0}:=H''(\cdot,\theta)$ we have
\begin{eqnarray}
i^f(B_{1})\notin  [i^f(B_{0}),i^f(B_{0})+\nu^f(B_{0}) ].\nn
\end{eqnarray}
Then (2.4) has a nontrivial periodic solution $x=x_0$.

Under the further assumption that

(iv) $\nu^f(B_{0})=0$ and $|i^f(B_{1})-i^f(B_{0})|\geq 2n$, (2.4)
has two nontrivial periodic solutions.

\vskip2mm
  {\bf Theorem 2.11}. Assume that

(i) $V''(t,x)$ is continuous and is bounded for $(t,x)\in
[0,1]\times \R^{2n}$, $V'(t,\theta)=\theta$;

(ii) there exist $B_{1}$, $B_{2}\in
L^{\infty}([0,1];GL_{s}(\R^{n}))$ satisfying
$i^s(B_{1})=i^s(B_{2})$, $\nu^s(B_{2})=0$ and
  $$
B_1(t)\leq V''(t,x)\leq B_2(t)\ \ \forall (t,x)\in [0,1]\times
\R^{n}\ \ with \ \ |x|\geq r>0;
 $$

(iii)with $B_{0}:=V''(\cdot,\theta)$ we have
\begin{eqnarray}
i^s(B_{1})\notin  [i^s(B_{0}),i^s(B_{0})+\nu^s(B_{0}) ].\nn
\end{eqnarray}
Then (2.9) has a nontrivial periodic solution.

Under the further assumption that

(iv) $\nu^s(B_{0})=0$ and $|i^s(B_{1})-i^s(B_{0})|\geq 2n$, (2.9)
has two nontrivial periodic solutions.

 {\bf Remark.} Theorems 2.6 and 2.9 were obtained already in [1] as special
 cases of a result concerning the first kind operator equation.
 However, the inequality (2.1) in [1] should be replaced by
   $$
 a(x,x)+\lambda_0 ||x||_X^2\geq c||x||_Z, \ \ \forall\ \ x\in Z
 $$
for some positive constants $\lambda_0$ and $c$.

 \setcounter{equation}{0}
\section{A new reduced functional}

In this section we will construct a new functional to investigate
(1.1). The method comes from Section 2.1 of Chapter IV in [4] and
[12]. Because every eigenvalue of $A$ is isolated, there exists
$\epsilon>0$ such that $A_{\epsilon}:=A+\epsilon I: D(A)\subset X\to
X$ is invertible and the inverse $A_{\epsilon}^{-1}: X\to X$
satisfies
    \be
||A_{\epsilon}^{-1}||\leq {1\over \epsilon}.\lb{3.1}
    \ee
Set $\Phi_{\epsilon}(x)=\Phi(x)+{1\over 2}\epsilon||x||^2$. Then
(1.1) is equivalent to the following equation
     \be
   A_{\epsilon}x-\Phi_{\epsilon}'(x)=\theta.
     \ee
Obviously $D(A_{\epsilon})=D(A), A_{\epsilon}:D(A)\subset X\to X$ is
selfadjoint and $\sigma(A_{\epsilon})=\sigma_d(A_{\epsilon})$. Let
$\{E'_{\lambda}\}$ be the spectral resolution of $A_{\epsilon}$.
There is an orthogonal decomposition:
     $$X=X^{+}\oplus X^{0}\oplus  X^{-}$$
where $X^{\ast}=P^{\ast}X$ for $\ast=+, 0, -$ and
$P^+=\int_{0}^{\infty} dE'_{\lambda}$, $P^{0}=\int_{-\beta}^{0}
dE'_{\lambda}$, $P^{-}=\int_{-\infty}^{-\beta} dE'_{\lambda}$ and
$\beta>0$. And from now on we always assume that $-\beta\in
\rho(A_{\epsilon})$. Let $x\in D(A_{\epsilon})$ be a solution of
(3.2). Set $x=x^{+}+x^{0}+x^{-},u =u^{+}+u^{0}+u^{-},
u^{\pm}=|A_{\epsilon}|^{\frac{1}{2}}x^{\pm},
u^0=|A_{\epsilon}|^{\frac{1}{2}}x^0$. Because
$A_{\epsilon}x=|A_{\epsilon}|(x^+-x^0-x^-)$,
$u=|A_{\epsilon}|^{\frac{1}{2}}x$ satisfies the following equation
   \be
u^{+}-u^{0}-u^{-}-|A_{\epsilon}|^{-\frac{1}{2}}\Phi_{\epsilon}'(|A_{\epsilon}|^{-\frac{1}{2}}u)
=\theta.
   \ee
Note that $V=D(|A|^{\frac{1}{2}})=D(|A_{\epsilon}|^{\frac{1}{2}})$.
Similar to page 189 in [4] by Chang, we define the functional as
follows
   \be
\varphi (u)=\frac{1}{2} \| u^{+} \|^{2}-\frac{1}{2}
\|u^{0}\|^{2}-\frac{1}{2} \| u^{-} \|^{2}-\Phi_{\epsilon}
(|A_{\epsilon}|^{-\frac{1}{2}}u), \ \ \forall \ \ u\in X. \lb{3.4}
   \ee
The Euler equation of this functional is (3.3). We only discuss
the case:
  \begin{eqnarray}
\sigma(A)\ \  \hbox {is unbounded from below}.
  \end{eqnarray}
After that we will find that if $\sigma(A)$ is bounded from below,
the things related are much simper. Since (3.5) holds,
dim$P^-X=+\infty$ and the Morse(negative) index at any critical
point is always infinite. In order to use Morse theory to
investigate (1.1) we need to obtain a reduced functional having a
finite Morse index at any critical point. To this end, we use the
method from [12]. Note that (3.3) is equivalent to the following
system:
   \begin{eqnarray}
&&u^{+}-u^{0}-(P^{+}+P^{0})|A_{\epsilon}|^{-\frac{1}{2}}\Phi_{\epsilon}'(|A_{\epsilon}|^{-\frac{1}{2}}u)=\theta\\
&&-u^{-}-P^{-}|A_{\epsilon}|^{-\frac{1}{2}}\Phi_{\epsilon}'(|A_{\epsilon}|^{-\frac{1}{2}}u)=\theta
   \end{eqnarray}
We will solve (3.7) for $u^{0}$, $u^{+}$ fixed. Denote the Frechet
derivative of $\Phi(x)$ with respect to $x$ in $V$ by $\Phi'_V(x)$.
Because $(\Phi'_V(u),v)_V=(\Phi'(u),v)$ for any $u,v\in V$, (3.7)
has an equivalent form
   \be
u^{-}=-P^{-}{(|A_{\epsilon}|^{-\frac{1}{2}})}^{\ast}{\Phi_{\epsilon}'}
_V(|A_{\epsilon}|^{-\frac{1}{2}}u).\lb{3.8}
  \ee
Set ${\cal N}(u^-)=-P^{-}(|A_{\epsilon}|^{-\frac{1}{2}})^{\ast}
{\Phi_{\epsilon}'}_V(|A_{\epsilon}|^{-\frac{1}{2}}(u^++u^0+u^-))$
for any $u^-\in X^-$. It suffices to prove that $\|{\cal
N}(u^-_1)-{\cal N}(u^-_2)\|\leq \alpha \|u^-_1-u^-_2\|$ for some
fixed $\alpha\in (0,1)$ and any $u^-_1,u^-_2\in X^-$. In fact, let
$\Phi''_V(x)$ denote the second Frechet derivative of $\Phi(x)$ with
respect to $x$ in $V$. It is also easy to check that
$(\Phi''_V(x)u,v)_V=(\Phi''(x)u,v)$ for any $x,u,v\in V$. Thus
   \begin{eqnarray}
&&{\cal N}(u^-_2)-{\cal
N}(u^-_{1})\nn\\
&&\ \
=P^{-}(|A_{\epsilon}|^{-\frac{1}{2}})^{\ast}{\Phi_{\epsilon}'}_{V}
(|A_{\epsilon}|^{-\frac{1}{2}}(u^{+}+u^{0}+u^-_1))-P^{-}(|A_{\epsilon}|^{-\frac{1}{2}})^{\ast}{\Phi_{\epsilon}'}_{V}
(|A_{\epsilon}|^{-\frac{1}{2}}( u^{+}+u^{0}+u^-_2))\nn\\
 &&\ \
  =P^{-}(|A_{\epsilon}|^{-\frac{1}{2}})^{\ast}\int_{0}^{1}{\Phi_{\epsilon}''}_{V}
                                                 (|A_{\epsilon}|^{-\frac{1}{2}}(u^{+}+u^{0}+\theta
                                                 u^-_1
+(1-\theta)u^-_2))d\theta|A_{\epsilon}|^{-\frac{1}{2}}(u^-_1-u^-_2)\nn\\
&&\ \
=P^{-}|A_{\epsilon}|^{-\frac{1}{2}}\int_{0}^{1}\Phi_{\epsilon}''
                                                 (|A_{\epsilon}|^{-\frac{1}{2}}(u^{+}+u^{0}+\theta
                                                 u^-_1
+(1-\theta)u^-_2))d\theta|A_{\epsilon}|^{-\frac{1}{2}}(u^-_1-u^-_2)\nn
   \end{eqnarray}
Let $\{e_{j} \}$ be the orthonormal basis of $X$ as in Section 1 and
$A_{\epsilon}e_{j}=\lambda^{'}_{j}e_{j}$ where
$\lambda_j^{'}=\lambda_j+\epsilon$ are all eigenvalues of
$A_{\epsilon}$ satisfying $\lambda^{'}_j\leq\lambda^{'}_{j+1} \ \
\forall j$ and $\lambda^{'}_j\to\pm\infty$ as $j\to\pm\infty$. Then
for any $x=\sum c_{j}e_{j}\in X$ we have
$P^{-}(|A_{\epsilon}|^{-\frac{1}{2}})x=\Sigma_{\lambda^{'}_{j<-\beta}}c_{j}(-\lambda^{'}_{j})^{-\frac{1}{2}}e_{j}$
and
   \begin{eqnarray}
\|P^{-}|A_{\epsilon}|^{\frac{1}{2}}x\|=(\Sigma_{\lambda^{'}_{j}<-\beta}
                  c_{j}^{2}(-\lambda^{'}_{j})^{-1})^{\frac{1}{2}}\leq
(\frac{1}{\beta}\Sigma
c_{j}^{2})^{\frac{1}{2}}=\frac{1}{\sqrt{\beta}}\|x\|.\nn
  \end{eqnarray}
Thus
  \be
\|P^{-}|A_{\epsilon}|^{-\frac{1}{2}}\|\leq \frac{1}{\sqrt{\beta}}.
  \ee
And by assumption(i) there exists $M>0$ such that
  \be
||\Phi''(x)||\leq M, \ \ ||\Phi_{\epsilon}''(x)||\leq M\ \ \forall \
\ x\in X.
  \ee
Hence $\|{\cal N}(u^-_2)-{\cal N}(u^-_{1})\|\leq \frac{M}{\beta}
\|u^-_2-u^-_{1}\|$. Let $\beta>0$ be large enough such that
$\frac{M}{\beta}<1$. Then ${\cal N}(u^-)=u^-$ and equivalently (3.7)
has a unique solution $u^-=u^{-}(u^+,u^0)\in C^{1}(X^{+}\oplus
X^{0}, X^{-})$.

Define
  \be
a(u^{+}+u^{0})=\varphi (u^{+}+u^{0}+u^{-}(u^{+},u^{0})).
  \ee
A critical point of $a(u^++u^0)$ corresponds to a solution of (3.2).
In fact,
  \begin{eqnarray}
a'(u^{+}+u^{0})&&=u^{+}-u^{0}-(u^{-'})^{\ast}
u^{-}-(P^{+}+P^{0}+(u^{-'})^{\ast}P^{-})(|A_{\epsilon}|^{-\frac{1}{2}})^{\ast}{\Phi_{\epsilon}'}_{V}
                                                                                (|A_{\epsilon}|^{-\frac{1}{2}}u)\nn\\
    &&=u^{+}-u^{0}-(P^{+}+P^{0})|A_{\epsilon}|^{-\frac{1}{2}}
                 {\Phi_{\epsilon}'}(|A_{\epsilon}|^{-\frac{1}{2}}u),
   \end{eqnarray}
where $u=u^{+}+u^{0}+u^{-}$ and $u^{-}$ satisfies (3.7). Hence,
$a'(u^{+}+u^{0})=0$ if and only if (3.6-3.7) hold and equivalently
(3.3) holds. Thus, we have the following proposition.

{\bf Proposition 3.1} Under assumptions (i-ii) of Theorem 1.1 the
functional $a(u^++u^0)$ defined in (3.11) belongs to $C^2(E)$, and
every critical point $u^++u^0$ corresponds to a solution
$x=|A_{\epsilon}|^{\frac{1}{2}}(u^++u^0+u^-(u^++u^0))$ of (3.2).

\vskip2mm In order to prove Theorem 1.1 we need to investigate the
Morse index of $a(u^++u^0)$ at a critical point. Let us calculate
$a''(u^++u^0)$ now. From (3.8) it follows that
\begin{eqnarray}
-P^{-}|A_{\epsilon}|^{-\frac{1}{2}}{\Phi_{\epsilon}''}
                       (|A_{\epsilon}|^{-\frac{1}{2}}u)|A_{\epsilon}|^{-\frac{1}{2}}(P^{+}+P^{0})
                       =(P^{-}+P^{-}|A_{\epsilon}|^{-\frac{1}{2}}{\Phi_{\epsilon}''}
                       (|A_{\epsilon}|^{-\frac{1}{2}}u)|A_{\epsilon}|^{-\frac{1}{2}}P^-)u^{-'}(u^{\ast}).\nn
  \end{eqnarray}
From (3.9-3.10) for $\beta>M$ the operator on the right side is
invertible and
  \begin{eqnarray}
u^{-'}(u^{\ast})
=&&-(P^{-}+P^{-}|A_{\epsilon}|^{-\frac{1}{2}}\Phi_{\epsilon}''(|A_{\epsilon}|^{-\frac{1}{2}}u)
                                 |A_{\epsilon}|^{-\frac{1}{2}}P^{-})^{-1}\nn\\
&&\ \ \ \
P^{-}|A_{\epsilon}|^{-\frac{1}{2}}
     \Phi''(|A_{\epsilon}|^{-\frac{1}{2}}u)|A_{\epsilon}|^{-\frac{1}{2}}(P^{+}+P^{0}).\nn
   \end{eqnarray}
Thus, (3.12) implies
  \begin{eqnarray}
a''(u^{+}+u^{0})&&=P^{+}-P^{0}-(P^{+}+P^{0})|A_{\epsilon}|^{-\frac{1}{2}}
                           \Phi''_{\epsilon}(|A_{\epsilon}|^{-\frac{1}{2}}u)|A_{\epsilon}|^{-\frac{1}{2}}(P^{+}+P^{0}+u^{-'}(u^{\ast}))\nn\\
   &&=P^{+}-P^{0}-(P^{+}+P^{0})|A_{\epsilon}|^{-\frac{1}{2}}\Phi_{\epsilon}''
                 (|A_{\epsilon}|^{-\frac{1}{2}}u)|A_{\epsilon}|^{-\frac{1}{2}}(P_{+}+P_{0})\nn\\
   &&+ (P^{+}+P^{0})|A_{\epsilon}|^{-\frac{1}{2}}\Phi_{\epsilon}''(|A_{\epsilon}|^{-\frac{1}{2}}u)
                            |A_{\epsilon}|^{-\frac{1}{2}}(P^{-}
\nn\\
   &&\ \ \ \ \  +P^{-}|A_{\epsilon}|^{-\frac{1}{2}}\Phi_{\epsilon}''
       (|A_{\epsilon}|^{-\frac{1}{2}}u)|A_{\epsilon}|^{-\frac{1}{2}}P^{-})^{-1}P^{-}
            |A_{\epsilon}|^{-\frac{1}{2}}\Phi_{\epsilon}''(|A_{\epsilon}|^{-\frac{1}{2}}u)
                                                         |A_{\epsilon}|^{-\frac{1}{2}}(P^{+}+P^{0}),\nn\\\lb{3.12}
   \end{eqnarray}
where $u=u^++u^0+u^-(u^+,u^0)$.

\vskip 4mm

In order to prove Theorem 1.1 we also need a lemma. Let $X$ be a
Hilbert space and $f\in C^2(X,\R)$. As in [4, Chapter 1] let
$K=\{x\in X|f'(x)=\theta\}, f_a=\{x\in X|f(x)\leq a\}$. If
$f'(x)=\theta, c=f(x)$ we say that $x$ is a critical point of $f$,
and $c$ is a critical value. $c\in \R$ is called a regular value of
$f$ if it is not a critical value. For any $x\in K$, $f''(x)\in
L_s(X)$ is selfadjoint. We call the dimension of the negative
subspace denoted by $m^-(f''(x))$ corresponding the the spectral
decomposing the Morse index of $x$, and $m^0(f''(x)):=\dim ker
f''(x)$ is called the Morse nullity of $x$. If $f''(x)$ has a
bounded inverse then $x$ is called non-degenerate. For any two
topological spaces $Y\subset X$ let $H_q(X,Y;\R)$ denote the $q$th
regular relative homology group. For an isolated critical point $x$,
the $q$th critical group is defined by $C_q(f,x)=H_q(f_c\cap
U,(f_c\setminus \{x\})\cap U);\R)$ for any neighborhood $U$ of $x$
with $U\cap K=\{x\}$ and $c=f(x_0)$. From [4, Chapter II Theorems
5.1 and 5.2] we have the following lemma.

{\bf Lemma 3.2}. Assume $f\in C^2(X,\R)$ satisfies the (PS)
condition, $f'(\theta)=\theta$, and there is a positive integer
$\gamma$  such that
$\gamma\notin[m^-(f''(\theta)),m^0(f''(\theta))+m^-(f''(\theta))]$
and $H_q(X,f_a;\R)=\delta_{q\gamma}\R$ for some regular value
$a<f(\theta)$. Then $f$ has a critical point $p_0\neq \theta$ with
$C_{\gamma}(f,p_0)\neq0$. Moreover, if $\theta$ is a
non-degenerate critical point, and
$m^0(f''(p_0))\leq|\gamma-m^-(f''(\theta))|$, then $f$ has another
critical point $p_1\neq p_0,\theta$.

\setcounter{equation}{0}
\section{Index theory for linear self-adjoint operator equations}

In this section we investigate the Morse index of $a''(u^\ast)$
obtained in the last section. This is a continuation to prepare for
the proof of Theorem 1.1. The method comes from [13-14]. For any
$B\in L_s(X)$, we set $B_{\epsilon}=B+\epsilon I$ where $\epsilon>0$
satisfies (3.1), and there exist large numbers $\beta>M>0$ such that
  \be
  \|B\|\leq M, \|B_{\epsilon}\|\leq M. \lb{4.1}
  \ee
Then
  \be
 \|P^-|A_{\epsilon}|^{-\frac{1}{2}}B_{\epsilon}|A_{\epsilon}|^{-\frac{1}{2}}P^-\|\leq
 {M\over\beta}\lb{4.2}
    \ee
and
$P^-+P^-|A_{\epsilon}|^{-\frac{1}{2}}B_{\epsilon}|A_{\epsilon}|^{-\frac{1}{2}}P^-:
X^-\to X^-$ is invertible. Motivated by (3.13) we consider the
following bilinear form
  \begin{eqnarray}
&&q_{A,\beta;B}(u^{\ast},v^{\ast})=\frac{1}{2}(u^{+},
v^{+})-\frac{1}{2}(u^{0},
v^{0})-\frac{1}{2}(|A_{\epsilon}|^{-\frac{1}{2}}B_{\epsilon}|A_{\epsilon}|^{-\frac{1}{2}}u^\ast,
v^\ast)\nn\\
       &&\ \ +\frac{1}{2}(|A_{\epsilon}|^{-\frac{1}{2}}B_{\epsilon}|A_{\epsilon}|^{-\frac{1}{2}}
       P^{-}(P^{-}
       +P^{-}|A|^{-\frac{1}{2}}B_{\epsilon}|A_{\epsilon}|^{-\frac{1}{2}}P^-)^{-1}P^{-}
       |A_{\epsilon}|^{-\frac{1}{2}}B_{\epsilon}|A_{\epsilon}|^{-\frac{1}{2}}u^\ast,v^\ast),
   \end{eqnarray}
where $u^{\ast}=u^++u^0, v^{\ast}=v^++v^0$ belong to
$E:=X^{+}\bigoplus X^{0}$. Define
    \begin{eqnarray}
    &&{\cal B}u^{\ast}=2P^{0}u^{\ast}+(P^++P^0)|A_{\epsilon}|^{-\frac{1}{2}}
             B_{\epsilon}|A_{\epsilon}|^{-\frac{1}{2}}u^{\ast}\nn\\
   &&-(P^++P^0)|A_{\epsilon}|^{-\frac{1}{2}}B_{\epsilon}|A_{\epsilon}|^{-\frac{1}{2}}
       (P^{-}+P^{-}|A_{\epsilon}|^{-\frac{1}{2}}B_{\epsilon}|A_{\epsilon}|^{-\frac{1}{2}}P^-)^{-1}P^{-}
       |A_{\epsilon}|^{-\frac{1}{2}}B_{\epsilon}|A_{\epsilon}|^{-\frac{1}{2}}u^{\ast}.\lb{4.1}
\end{eqnarray}
Then
$q_{A,\beta;B}(v^{\ast},u^{\ast})=\frac{1}{2}[(v^{\ast},u^{\ast})-({\cal
B}v^{\ast},u^{\ast})] $ for any $v^{\ast}, u^{\ast}\in E$. And
${\cal B}:E \rightarrow E$ is self-adjoint and compact. By the
spectral theory there is a basis $\{e_{j}\}$ of $E$ and a sequence
$\mu_{j}\rightarrow 0$ in $\textbf{R}$ such that:
  \begin{eqnarray}
{\cal B}e_{j}=\mu_{j}e_{j}; \hbox{ } (e_{j},e_{i})=\delta_{ij} .
  \end{eqnarray}
For any $u^{\ast}\in E$, which can be expressed as
$u^{\ast}=\sum^{\infty}_{j=1}c_{j} e_{j}$
    \begin{eqnarray}
q_{A,\beta;B}(u^{\ast},u^{\ast}) =\frac{1}{2} \sum _{j=1}^{\infty}
               (1-\mu_{j})c_{j}^{2}.\nn
  \end{eqnarray}
Define
  \begin{eqnarray}
E^{-}_{A,\beta}(B):&=&\{\sum_{j=1}^{\infty}c_{j}e_{j}\mid c_{j}=0 \hbox{ if } 1-\mu_{j}\geq 0  \},\nn\\
E^{0}_{A,\beta}(B):&=&\{\sum_{j=1}^{\infty}c_{j}e_{j}\mid c_{j}=0 \hbox{ if } 1-\mu_{j}\neq 0  \},\nn\\
E^{+}_{A,\beta}(B):&=&\{\sum_{j=1}^{\infty}c_{j}e_{j}\mid c_{j}=0
\hbox{ if } 1-\mu_{j}\leq 0  \}.\nn
   \end{eqnarray}
Obviously, $E^{-}_{A,\beta}(B)$, $E^{0}_{A,\beta}(B)$ and
$E^{+}_{A,\beta}(B)$ are $q_{A,\beta;B}$-orthogonal and
$E^{-}_{A,\beta}(B)\oplus E^{0}_{A,\beta}(B)\oplus
E^{+}_{A,\beta}(B)=E$. Since $\mu_{j}\rightarrow 0$ as
$j\rightarrow \infty$, $E^{-}_{A,\beta}(B)$ and
$E^{0}_{A,\beta}(B)$ are two finite dimensional subspaces. \vskip
2mm {\bf Definition 4.1}. We define
$$\nu_{A,\beta}(B):=\hbox{dim } E^{0}_{A,\beta}(B), i_{A,\beta}(B):=\hbox{dim }
E^{-}_{A,\beta}(B).$$

\vskip4mm

 {\bf Proposition 4.2}. We have the following results:

(i)$\nu_{A,\beta}(B)=dim(\hbox{ker}(A-B))$;

(ii)$i_{A,\beta}(B)$ is the Morse index of $q_{A,\beta;B}$;

 (iii) For any $B_0, B_1\in L_s(X)$ satisfying
$B_0\leq B_1$ and $B_0<B_1$ with respect to $ker(A-B_{\lambda})$
$\forall \ \lambda\in [0,1)$ where $B_{\lambda}=
B_0+\lambda(B_1-B_0)$ if the subspace is not trivial, we have
   $$
i_{A,\beta}(B_1)-i_{A,\beta}(B_0)=\Sigma_{\lambda\in
[0,1)}\nu_{A,\beta}(B_{\lambda}).
  $$

{\bf Proof}. (i) Fix any $u^{\ast}\in E^{0}_{A,\beta}(B)$; by
definition $q_{A,\beta;B}(u^{\ast},v^{\ast})=0\ \ \forall \
v^{\ast}\in E$. It follows
  \begin{eqnarray}
(P^{+}-P^{0})u^{\ast}-(P^{+}+P^{0})|A_{\epsilon}|^{-\frac{1}{2}}B_{\epsilon}|A_{\epsilon}|^{-\frac{1}{2}}u^{\ast}
-(P^{+}+P^{0})|A_{\epsilon}|^{-\frac{1}{2}}B_{\epsilon}|A_{\epsilon}|^{-\frac{1}{2}}u^{-}=0,\lb{4.5}
  \end{eqnarray}
where
  \begin{eqnarray}
u^{-}=-(P^{-}+P^{-}|A_{\epsilon}|^{-\frac{1}{2}}B_{\epsilon}
        |A_{\epsilon}|^{-\frac{1}{2}}P^{-})^{-1}P^{-}|A_{\epsilon}|^{-\frac{1}{2}}
             B_{\epsilon}|A_{\epsilon}|^{-\frac{1}{2}}u^{\ast},\nn
\end{eqnarray}
from which it follows
\begin{eqnarray}
u^{-}
=-P^{-}|A_{\epsilon}|^{-\frac{1}{2}}B_{\epsilon}|A_{\epsilon}|^{-\frac{1}{2}}(u^{\ast}+u^{-}).\nn
 \end{eqnarray}
Let $x^{-}=|A_{\epsilon}|^{-\frac{1}{2}}u^{-},
x^{+}+x^{0}=|A_{\epsilon}|^{-\frac{1}{2}}u^{\ast}$, and
$x=x^++x^0+x^-$. We obtain
  \begin{eqnarray}
A_{\epsilon}x^{-}-P^{-}B_{\epsilon}x=0.
   \end{eqnarray}
From (4.6) we obtain
\begin{eqnarray}
(P^{+}-P^{0})u^{\ast}-(P^{+}+P^{0})|A_{\epsilon}|^{-\frac{1}{2}}B_{\epsilon}(|A_{\epsilon}|^{-\frac{1}{2}}u^{\ast}
      +|A_{\epsilon}|^{-\frac{1}{2}}u^{-})=0.\nn
\end{eqnarray}
That is
\begin{eqnarray}
A_{\epsilon}(x^{+}+x^{0})-(P^{+}+P^{0})B_{\epsilon}x=0.\lb{4.7}
\end{eqnarray}
Combining (4.7) and (4.8) implies
    $$Ax-Bx=0.$$
Thus, from the one to one correspondence $u^{\ast}\mapsto
x=|A_{\epsilon}|^{-\frac{1}{2}}(u^{\ast}+u^{-})$ it follows that
$E^{0}_{A,\beta}(B)\cong \hbox{ker}(A-B)$.

(ii) Assume $X_{1}$ is a subspace of $E$ such that $q_{A,\beta;B}$
is negative definite on $X_{1}$ with dim$X_{1}=k$. Let $\{ x_{j}
\}_{1}^{k}$ be linear independent in $X_{1}$. We have the
decomposition $x_{j}=x_{j}^{-}+x_{j}^{\ast}$ with $x_{j}^{-}\in
E^{-}_{A,\beta}(B)$ and $x_{j}^{\ast} \in E^{0}_{A,\beta}(B)\oplus
E^{+}_{A,\beta}(B)$. If there exist not all zero numbers
$\alpha_{i}\in \textbf{R}$ such that
$\sum_{i=1}^{k}\alpha_{i}x_{i}^{-}=\theta$. On the one hand,
$x:=\sum_{i=1}^{k}\alpha_{i}x_{i}\in X_{1}\setminus\{\theta\}$ and
$q_{A,\beta;B}(x,x)< 0$; on the other hand,
$x=\sum_{i=1}^{k}\alpha_{i}x_{i}^{\ast}\in E^{0}_{A,\beta}(B)\oplus
E^{+}_{A,\beta}(B)$, and $q_{A,\beta;B}(x,x)\geq 0$. So
$\{x_{i}^{-}\}_{i=1}^{k}$ is linear independent and
$i_{A,\beta}(B)\geq k=$dim$X_{1}$.

(iii) From Theorem 1.5(ii) and (1.4), $
i_{A}(B_1)-i_{A}(B_0)=\Sigma_{\lambda\in
[0,1)}\nu_{A}(B_{\lambda})$, and there are at most finite numbers
$\lambda\in[0,1)$ such that $ker(A-B_{\lambda})\neq \{\theta\}$. Set
$i(\lambda)=i_{A,\beta}(B_{\lambda}),\
\nu(\lambda)=\nu_{A,\beta}(B_{\lambda})$ for $\lambda\in [0,1)$. Let
${\cal B}(\lambda)$ be the operator similar to ${\cal B}$ in (4.4)
only with $B_{\epsilon}$ replaced with $B_{\lambda}+\epsilon I$. We
have the following lemma. \vskip2mm {\bf Lemma 4.3}. Let $M>0$ be
large enough such that $\| B_{0}\|\leq M$, $\| B_{1}\|\leq M$,
$\epsilon <M$ and $4M+\epsilon<\beta$. Then we have the following
expression
  \begin{eqnarray}
{\cal B}(\lambda)={\cal B}(\lambda_{0})+(\lambda-\lambda_{0}){\cal
B}_{1}(\lambda_0)+\cdots +(\lambda-\lambda_{0})^{k}{\cal
B}_{k}(\lambda_0)+\cdots, \nn
   \end{eqnarray}
where ${\cal B}_{k}:E\rightarrow E$ is selfadjoint and compact and
satisfies $$\|{\cal B}_{k}(\lambda_0) \|\leq \frac{8M}{\epsilon},\
\ k=2, 3, 4\cdots,$$
 and
   \begin{eqnarray}
&&{\cal B}_{1}(\lambda_{0})=(P^{+}+P^{0})(I+S^{T})|A_{\epsilon}|^{-\frac{1}{2}}(B_{1}-B_{0})|A_{\epsilon}|^{-\frac{1}{2}}(I+S)(P^{+}+P^{0}),\nn\\
&&S=-P^{-}(P^{-}+P^{-}|A_{\epsilon}|^{-\frac{1}{2}}B_{\lambda_{0}}|A_{\epsilon}|^{-\frac{1}{2}}P^{-})^{-1}
P^{-}|A_{\epsilon}|^{-\frac{1}{2}}B_{\lambda_{0}}|A_{\epsilon}|^{-\frac{1}{2}}.\nn
   \end{eqnarray}
{\bf Proof}. We have
  \begin{eqnarray}
&&P^{-}+P^{-}|A_{\epsilon}|^{-\frac{1}{2}}(B_{\lambda})_{\epsilon}|A_{\epsilon}|^{-\frac{1}{2}}P^{-}\nn\\
&&=(P^{-}+P^{-}|A_{\epsilon}|^{-\frac{1}{2}}(B_{\lambda_0})_{\epsilon}|A_{\epsilon}|^{-\frac{1}{2}}P^{-})
 (P^{-}\nn\\
 &&\ \ \ \ \ +(\lambda-\lambda_0)(P^{-}+P^{-}|A_{\epsilon}|^{-\frac{1}{2}}(B_{\lambda_0})_{\epsilon}
     |A_{\epsilon}|^{-\frac{1}{2}}P^{-})^{-1}
P^{-}|A_{\epsilon}|^{-\frac{1}{2}}(B_{1}-B_0)|A_{\epsilon}|^{-\frac{1}{2}}P^{-}).\nn\\
&&:=Q_{1}(P^-+(\lambda-\lambda_0)Q_{1}^{-1}Q_{2}),\nn
   \end{eqnarray}
where $Q_{1}=P^{-}+
P^{-}|A_{\epsilon}|^{-\frac{1}{2}}(B_{\lambda_0})_{\epsilon}
     |A_{\epsilon}|^{-\frac{1}{2}}P^{-}$, $Q_{2}=
P^{-}|A_{\epsilon}|^{-\frac{1}{2}}(B_1-B_0)
     |A_{\epsilon}|^{-\frac{1}{2}}P^{-}$, $\|Q_1^{-1}\|\leq \frac{1}{1-\frac{M+\epsilon}{\beta}}$
 ,  $\|Q_2\|\leq \frac{2M}{\beta}$.
  \begin{eqnarray}
&&(P^{-}+P^{-}|A_{\epsilon}|^{-\frac{1}{2}}(B_{\lambda})_{\epsilon}|A_{\epsilon}|^{-\frac{1}{2}}P^{-})^{-1}\nn\\
&&=(P^-+(\lambda-\lambda_0)Q_{1}^{-1}Q_{2})^{-1}Q_{1}^{-1}\nn\\
&&=Q_{1}^{-1}-(\lambda-\lambda_0)Q_{1}^{-1}Q_{2}Q_{1}^{-1}+(\lambda-\lambda_0)^{2}(Q_{1}^{-1}Q_{2})^{2}Q_{1}^{-1}+\cdots\nn
   \end{eqnarray}
Thus, the third term in ${\cal B}(\lambda)$ is
\begin{eqnarray}
&&-Q_{3}Q_{1}^{-1}Q_{3}^{T}+(\lambda-\lambda_{0})[Q_{3}Q_{1}^{-1}Q_{2}Q_{1}^{-1}Q_{3}^{T}
-Q_{3}Q_{1}^{-1}Q_{4}^{T}-Q_{4}Q_{1}^{-1}Q_{3}^{T}]+\cdots \nn\\
&&+(\lambda-\lambda_{0})^{k}[(-1)^{k+1}Q_{3}(Q_{1}^{-1}Q_{2})^{k}Q_{1}^{-1}Q_{3}^{T}
+(-1)^{k}Q_{4}(Q_{1}^{-1}Q_{2})^{k-1}Q_{1}^{-1}Q_{3}^{T}\nn\\
&&+(-1)^{k}Q_{3}(Q_{1}^{-1}Q_{2})^{k-1}Q_{1}^{-1}Q_{4}^{T}+(-1)^{k-1}Q_{4}(Q_{1}^{-1}Q_{2})^{k-2}Q_{1}^{-1}Q_{4}^{T}]+\cdots,\nn
\end{eqnarray}
where
$Q_{3}=(P^++P^0)|A_{\epsilon}|^{-\frac{1}{2}}(B_{\lambda_{0}})_{\epsilon}|A_{\epsilon}|^{-\frac{1}{2}}P^{-}$,
$Q_{4}=(P^++P^0)|A_{\epsilon}|^{-\frac{1}{2}}(B_{1}-B_{0})|A_{\epsilon}|^{-\frac{1}{2}}P^{-}$.
 Then
$\|Q_{3}\|\leq\frac{2M}{\sqrt{\epsilon\beta}}$,
$\|Q_{4}\|\leq\frac{2M}{\sqrt{\epsilon\beta}}$,
$\|Q_{1}^{-1}\|\leq\frac{1}{1-\frac{M+\epsilon}{\beta}}$ and
$\|Q_{2}\|\leq\frac{2M}{\beta}$. Thus, the results follow.

\vskip2mm Now we give the proof in three steps.

{\bf Step 1}. If $\nu(\lambda)=0$, then $i(\lambda)$ is continuous.
Suppose that $\nu(\lambda_0)=0$ for some $\lambda_0\in(0,1)$. Let
$S_1$ be the unit ball of $E^-_{A,\beta}(B_{\lambda_0})$. Because
the subspace is finite dimensional, $S_1$ is compact and
$q_{A,\beta;B_{\lambda}}(u,u)$ is continuous with respect to
$(\lambda,u)\in [0,1]\times S_1$. Thus, for $\lambda$ close enough
to $\lambda_0$, $q_{A,\beta;B_{\lambda}}$ is negative definite on
$S_{1}$ and hence on $E^-_{A,\beta}(B_{\lambda_0})$. By (ii),
$i(\lambda)\geq i(\lambda_0)$. In the following we prove the inverse
inequality. If $i(\lambda_l)>i(\lambda_0):=k$ for $\lambda_l\to
\lambda_0$,  then similar to (4.5) we have
   \begin{eqnarray}
{\cal B}(\lambda_l)e_{l,j}=\mu_{l,j}e_{l,j}; \hbox{ }
(e_{l,j},e_{l,i})=\delta_{ij}  \hbox{ for any }1\leq i,j  \leq
k+1,
   \end{eqnarray}
where  $\hbox{span}\{ e_{l,j}\}\subseteq
E_{A,\beta}^{-}(B_{\lambda_l})$. By definition $\mu_{l,j}=({\cal
B}(\lambda_l)e_{l,j},e_{l,j})$ is bounded in $\textbf{R}$ for
$j=1,\cdots k+1$, and $\lambda_{l}\in[0,1]$. So we can assume
 $e_{l,j}\rightharpoonup e_{j}$,
$\mu_{l,j}\rightarrow \mu_{j}$ and ${\cal
B}(\lambda_l)e_{l,j}\rightarrow {\cal B}(\lambda_0)e_{j}$ by gong
to subsequences if necessary. Taking the limit in (4.9) we obtain
  \begin{eqnarray}
{\cal B}(\lambda_0)e_{j}=\mu_{j}e_{j} .
  \end{eqnarray}
By definition for $j=1,\cdots k+1$, $1+\mu_{l,j}<0$ and
$\{\frac{1}{\mu_{l,j}} \}$ is bounded in $\textbf{R}$. So
  \begin{eqnarray}
e_{l,j}=\frac{1}{\mu_{l,j}} {\cal B}(\lambda)e_{l,j}\rightarrow
\frac{1}{\mu_{j}} {\cal B}(\lambda_{0})e_{j}=e_{j}  \nn
  \end{eqnarray}
in $E$, and  $(e_{i}, e_{j})=\delta_{ij}$ $1\leq i,j\leq k+1$. It
follows that $\{e_j\}_{1}^{k+1}$ is independent. And for every
$u=\sum_{i=1}^{k+1}c_{j}e_{j}$, since
$\sum_{i=1}^{k+1}c_{j}e_{l,j}\rightarrow u$ in $E$ and
 \begin{eqnarray}
q_{A,\beta;B_{\lambda_l}}(\sum_{i=1}^{k+1}c_{j}e_{l,j},\sum_{i=1}^{k+1}c_{j}e_{l,j})<0,\nn
  \end{eqnarray}
taking the limit as $l\rightarrow \infty$ we have
  \begin{eqnarray}
q_{A,\beta;B_{\lambda_0}}(u,u)\leq 0.\nn
 \end{eqnarray}
This means that $i(\lambda_0)\geq k+1$, a contradiction.

{\bf Step 2}. If $ \nu (\lambda_0)=0$ does not hold for
$\lambda_0\in [0,1)$, then
$i(\lambda_0+0)=i(\lambda_0)+\nu(\lambda_0)$. From the argument
above it follows that $i(\lambda_0+0)\leq
i(\lambda_0)+\nu(\lambda_0)$. Thus, by (ii) we need only to prove
that $q_{A,\beta;B_{\lambda}}(u,u)\leq 0$ for every $u\in
E^{-}_{A,\beta}(B_{\lambda_0})\oplus
E^{0}_{A,\beta}(B_{\lambda_0})\backslash \{\theta\}$ with
$\lambda-\lambda_0>0$ small. Let $S=\{u\in
E^{-}_{A,\beta}(B_{\lambda_0})\oplus
E^{0}_{A,\beta}(B_{\lambda_0}) | \| u \| =1 \}$. Because $S$ is
compact, we need only to prove that $\forall \ \ u^{\ast}\in S$
there exists $\delta>0$ that for any
$\delta>\lambda-\lambda_{0}>0$,
$q_{A,\beta;B_{\lambda}}(v^{\ast},v^{\ast})<0$ if $v^{\ast}$ is
close to $u^{\ast}$.

 {\bf Case 1}. Assume $u^{\ast}=u^{-}+u^{0}\in S$, $u^{-}\neq \theta$.
It follows
 \begin{eqnarray}
q_{A,\beta;B_{\lambda}}(v^{\ast},v^{\ast})=
q_{A,\beta;B_{\lambda_{0}}}(v^{\ast},v^{\ast})+o(1)\nn
  \end{eqnarray}
as $\lambda\rightarrow \lambda_0^{+}$.

Because $q_{A,\beta;B_{\lambda_{0}}}(u^{-},u^{-})<0$, there exists a
neighborhood $U$ of $u^{\ast}$  in $S$ such that $\forall \ \
v^{\ast}\in U$
\begin{eqnarray}
q_{A,\beta;B_{\lambda_{0}}}(v^{\ast},v^{\ast})<\frac{1}{2}q_{A,\beta;B_{\lambda_{0}}}(u^{-},u^{-}).\nn
\end{eqnarray}
The results follows.

 {\bf Case 2}. Assume $u^{\ast}\in S\bigcap E^{0}_{A,\beta}(B_{\lambda_0})$.

By the proof of (i) $({\cal
B}_{1}(\lambda_{0})u^{\ast},u^{\ast})>0$. There exists a
neighborhood $U$ in $S$ such that $\forall v^{\ast}\in U$, $({\cal
B}_{1}(\lambda_{0})v^{\ast},v^{\ast})>\frac{1}{2}({\cal
B}_{1}(\lambda_{0})u^{\ast},u^{\ast})$. It follows from Lemma 4.3
that
\begin{eqnarray}
q_{A,\beta;B_{\lambda}}(v^{\ast},v^{\ast})&&\leq
q_{A,\beta;B_{\lambda_{0}}}(v^{\ast},v^{\ast})-\frac{1}{2}(\lambda-\lambda_{0})({\cal B}_{1}(\lambda_{0})v^{\ast},v^{\ast})+o(\lambda-\lambda_{0})\nn\\
    &&\leq
    -\frac{1}{4}(\lambda-\lambda_{0})({\cal B}_{1}(\lambda_{0})u^{\ast},u^{\ast})+o(\lambda-\lambda_{0})\nn\\
    &&<0, \    \ \hbox{as}\lambda\rightarrow\lambda_{0}^+.\nn
  \end{eqnarray}

{\bf Step 3}. If $\nu(\lambda_0)\neq 0$  for same
$\lambda_0\in(0,1)$, then $i(\lambda_0-0)=i(\lambda_0)$. In fact,
suppose $i(\lambda_0-0)>i(\lambda_0)$. As in Step 1, there exists
$\lambda_l\rightarrow\lambda_0^-$ and $\mu_{l}\in\sigma({\cal
B}(\lambda_l))$ satisfying $\mu_{l}\rightarrow  1^{+}$. However,
Lemma 4.3 tells us
   $$
{\cal B}(\lambda)={\cal B}(\lambda_0)+(\lambda-\lambda_0){\cal
B}_1(\lambda_0)+o(\lambda-\lambda_0).
   $$
By Rellich's theory(Theorem 3.9 and Remark 3.11 in pages 392-393
of [15]), its eigenvalue and associated eigenvector are
holomorphic:
  \begin{eqnarray}
&&\mu(\lambda)=1+(\lambda-\lambda_0)\mu_1+o(\lambda-\lambda_0),\nn\\
&&u(\lambda)=u_0+(\lambda-\lambda_0)u_1++o(\lambda-\lambda_0),
B(\lambda_{0})u_{0}=u_{0}\neq \theta. \nn
  \end{eqnarray}
Thus, $\mu_1=({\cal B}_1(\lambda_0)u_0,u_0)>0$ and
$\mu_l=\mu(\lambda_l)=1+(\lambda_{l}-\lambda_{0})\mu_{1}+o(\lambda_{l}-\lambda_{0})<1$
as $\lambda_{l}\rightarrow\lambda_{0}^{-}$. This is a
contradiction. \hfill\hb

Note that $q_{A,\beta;B}(u,u)$ is not monotone with respect to $B$
because of the last term. This is not convenient to the proof of
Theorem 1.1. However, the inequality (3.9) implies that the last
term is much smaller than the third term if $\beta>0$ is large
enough. Thus, we can use the sum of the first three terms denoted
by $\bar{q}_{A,\beta;B}(u,u)$ instead, and we have the following
proposition.

\vskip 4mm {\bf Proposition 4.3} (i) Assume that $B\in L_s(X)$
satisfies $\nu_{A,\beta}(B)=0$. Then for $\beta>0$ large enough,
$({\bar q}_{A,\beta;B}(u,u))^{1\over 2}$ and $(-{\bar
q}_{A,\beta;B}(u,u))^{1\over 2}$ are equivalent norms on
$E^+_{A,\beta}(B)$ and $E^-_{A,\beta}(B)$ respectively.

(ii) Assume that $B_1,B_2\in L_s(X)$ satisfy $B_1<B_2,
i_{A,\beta}(B_1)=i_{A,\beta}(B_2)$ and
$\nu_{A,\beta}(B_1)=i_{A,\beta}(B_2)=0$. Then, for $\beta>0$ large
enough $E=E^+_{A,\beta}(B_2)\bigoplus E^-_{A,\beta}(B_1)$.

{\bf Proof}. Recall that $P^{+}=\int_{0}^{+\infty} dE^{'}_{\lambda}
$, $P^{0}=P^{0}_{\beta}=\int_{-\beta}^{0} dE^{'}_{\lambda} $,
 $P^{-}=P^{-}_{\beta}=\int_{-\infty}^{-\beta}  dE^{'}_{\lambda}
X$ and the operator ${\cal B}$ defined in (4.4) depends on $\beta$,
so we denote it by ${\cal B}_{\beta}$ now.

 (i) We need only to prove that there exists
$\delta>0$ such that
    \begin{eqnarray}
    (1-\delta,1+\delta)\bigcap \sigma({\cal B}_\beta)=\emptyset
  \end{eqnarray}
if $\beta>\beta_0$ for some $\beta>0$ large enough. Otherwise,
there exist $\beta_{k}\rightarrow +\infty, \mu_k\rightarrow 1$ and
unit vector $e_{k}^{\ast}\in E$ satisfying
 \begin{eqnarray}
{\cal B}_{\beta_{k}}e_{k}^{\ast}=\mu_k e_{k}^{\ast}.\nn
  \end{eqnarray}
By definition we have
 \begin{eqnarray}
(\mu_{k}-2)P^{0}_{k}e_{k}^{\ast}+\mu_{k}P^{+}e_{k}^{\ast}=
(P^{+}+P^{0}_{k})|A_{\epsilon}|^{-\frac{1}{2}}B_{\epsilon}|A_{\epsilon}|^{-\frac{1}{2}}e_{k}^{\ast}
-(P^{+}+P^{0}_{k})|A_{\epsilon}|^{-\frac{1}{2}}B_{\epsilon}|A_{\epsilon}|^{-\frac{1}{2}}e_{k}^{-},
  \end{eqnarray}
where $e_{k}^{-}=(
P^{-}_{k}+P^{-}_{k}|A_{\epsilon}|^{-\frac{1}{2}}B_{\epsilon}|A_{\epsilon}|^{-\frac{1}{2}}P^{-}_{k}
)^{-1}P^{-}_{k}|A_{\epsilon}|^{-\frac{1}{2}}B_{\epsilon}|A_{\epsilon}|^{-\frac{1}{2}}e_{k}^{\ast}$,
 $P^{0}_{k}=P_{\beta_{k}}^{0}$, $P^{-}_{k}=P_{\beta_{k}}^{-}$ and
 ${\cal B}_{k}={\cal B}_{\beta_{k}}$.
 Assume $e_{k}^{\ast}\rightharpoonup  e=e^{+}+e^{0} $, then
 $e_{k}^{+}\rightharpoonup e^+$, $e_{k}^{0}\rightharpoonup e^0$
 and
$|A_{\epsilon}|^{-\frac{1}{2}}e_{k}^{\ast}\rightarrow
|A_{\epsilon}|^{-\frac{1}{2}}e$. By (4.11), $e_{k}^{+}\rightarrow
e^+=P^{+}|A_{\epsilon}|^{-\frac{1}{2}}B_{\epsilon}|A_{\epsilon}|^{-\frac{1}{2}}e$,
 $e_{k}^{0}\rightarrow e^0=(I-P^{+})|A_{\epsilon}|^{-\frac{1}{2}}B_{\epsilon}|A_{\epsilon}|^{-\frac{1}{2}}e$.
Thus  $e_{k}^{\ast}\rightarrow e$ and $\| e\|= \|e_{k}^{\ast}\|=1$.
On the other hand, it follows that
$e^++e^{0}=|A_{\epsilon}|^{-\frac{1}{2}}B_{\epsilon}|A_{\epsilon}|^{-\frac{1}{2}}e$,
and $A_{\epsilon}x=B_{\epsilon}x$ with
$x=|A_{\epsilon}|^{-\frac{1}{2}}e\neq 0$, a contradiction.

By definition, for $\beta\geq\beta_{0}$ and
$u=\sum_{1-\mu_j>0}c_{j}e_{j}\in E^+_{A,\beta}(B)$ it follows
   \begin{eqnarray}
 \hbox{sup} \{1-\mu_{j}| 1-\mu_{j}>0  \} ||u||^2 \geq
q_{A,\beta;B}(u,u))\geq \delta ||u||^{2}.
   \end{eqnarray}
And by (4.1), (4.2), (4.3),
  \begin{eqnarray}
|q_{A,\beta;B}(u,u)-\bar{q}_{A,\beta;B}(u,u)
|\leq\frac{M^2}{\epsilon(\beta-M)}\|u\|^{2}.
  \end{eqnarray}
Thus, for $\beta>0$ large enough,
$(\bar{q}_{A,\beta;B}(u,u))^{1\over2}$ is an equivalent norm in
$E^+_{A,\beta}(B)$.

(ii) Assume $u$ belongs both $E^+_{A,\beta}(B_2)$ and
$E^-_{A,\beta}(B_1)$. By definition and the result in (i),
 $$
 0\geq {\bar q}_{A,\beta;B_1}(u,u)\geq {\bar q}_{A,\beta;B_2}(u,u)\geq 0.
  $$
Thus, $u=\theta$. Now we need only to prove that
$E=E^+_{A,\beta}(B_2)+E^-_{A,\beta}(B_1)$. In fact, let $\{ e_{j}
\}_{1}^{k}$ be a basis of $E^-_{A,\beta}(B_1)$. We have the
decomposition $e_{j}=e_{j}^{-}+e_{j}^{+}$ with $e_{j}^{-}\in
E^{-}_{A,\beta}(B_2)$ and $e_{j}^{+} \in E^{+}_{A,\beta}(B_2)$. If
there exist not all zero numbers $c_{j}\in \textbf{R}$ such that
$\sum_{j=1}^{k}c_{j}e_{j}^{-}=\theta$. On the one hand,
$x:=\sum_{i=1}^{k}c_{j}e_{j}\in E^{-}_{A,\beta}(B_1)\setminus
\{\theta\}$ and ${\bar q}_{A,\beta;B_1}(x,x)< 0$ if $\beta>0$ is
large enough; on the other hand, $x=\sum_{j=1}^{k}c_{j}e_{j}^{+}\in
E^{+}_{A,\beta}(B_2)$, and ${\bar q}_{A,\beta;B_2}(x,x)\geq 0$ if
$\beta>0$ is large enough. This is a contradiction. So
$\{e_{j}^{-}\}_{j=1}^{k}$ is linear independent. For any $u\in E$,
$u=u^-+u^+$ with $u^-\in E^{-}_{A,\beta}(B_2),u^+\in
E^{+}_{A,\beta}(B_2)$. There exist $\{c_{j} \}_{j=1}^{k} \subset\psi
\textbf{R}$ such that $u^-=\sum_{j=1}^{k}c_{j}e_{j}^{-}$. Thus,
$u=\sum_{j=1}^{k}c_{j}e_{j}+(u^+-\sum_{j=1}^{k}c_{j}e_{j}^+)$.\hfill\hb

\setcounter{equation}{0}
\section{Proof of Theorem 1.1}

 In this section we prove Theorem 1.1. We have the  following two propositions.

{\bf Proposition 5.1}. Under assumptions (i-ii) of Theorem 1.1
$a(u^{\ast})$ satisfies the (PS) condition.

{\bf Proposition 5.2}. For $c>0$(such that $-c<f(0)$)  large enough,
we have
  \begin{eqnarray}
H_{q}(E;a_{-c}, \textbf{R}) \cong \delta_{q\gamma}\R, \hbox{
}q=0,1,2,\cdots,
   \end{eqnarray}
where $\gamma=i_{A,\beta}(B_{1})$.

{\bf Proof of Proposition 5.1}. Assume that $\{ u_{n}^{+}+ u_{n}^{0}
\}$ is a
 sequence in $E$ such that  $a'(u_{n}^{+}+ u_{n}^{0} )\rightarrow 0$ in $E$.
From (3.12),
   \begin{eqnarray}
&&a'(u^{+}_{n}+ u^{0}_{n} )=u^{+}_{n}-
u^{0}_{n}-(P^{+}+P^{0})|A_{\epsilon}|^{-\frac{1}{2}}\Phi_{\epsilon}
' (|A_{\epsilon}|^{-\frac{1}{2}}u_{n})\nn\\
&&u^{-}_{n}=-P^{-}|A_{\epsilon}|^{-\frac{1}{2}}\Phi_{\epsilon} '
(|A_{\epsilon}|^{-\frac{1}{2}}u_{n}),\ \ u_n=u_n^++u_n^0+u_n^-\nn.
   \end{eqnarray}
We claim that $\{ u_{n}^{+}+ u_{n}^{0} \}$ is bounded. If the case
is not true, then $||u_n||\geq\|  u_{n}^{+}+ u_{n}^{0}
\|\rightarrow\infty$. From assumption (ii) it follows
$\Phi_{\epsilon} '
(|A_{\epsilon}|^{-\frac{1}{2}}u_n)=B_{\epsilon}(|A_{\epsilon}|^{-\frac{1}{2}}u_n)|A_{\epsilon}|^{-\frac{1}{2}}u_n+C_n$
satisfying $C_n\in X$ is bounded and
  $$
B_1\leq B(|A_{\epsilon}|^{-\frac{1}{2}}u_n)\leq B_2.
  $$
Let
$y_{n}=u_{n}/
\|u_{n}\|$. Then
  \begin{eqnarray}
y^{+}_{n}-y^{0}_{n}-y^-_n-|A_{\epsilon}|^{-\frac{1}{2}}B_{\epsilon}(|A_{\epsilon}|^{-\frac{1}{2}}u_{n})|A|^{-\frac{1}{2}}
y_{n}- \|u_{n}\|^{-1}|A_{\epsilon}|^{-\frac{1}{2}}C_n\rightarrow 0.
    \end{eqnarray}
Because $||y_n||=1$, $B(|A_{\epsilon}|^{-\frac{1}{2}}u_n)\in
L_s(X)$ is bounded and $X$ is separable, we can assume as in [1,
page 81] that $y_{n}\rightharpoonup y$ in $X$ and
$B(|A_{\epsilon}|^{-\frac{1}{2}}u_n)x\rightharpoonup Bx$ for any
$x\in X$ and some $B\in L_s(X)$ such that $B_1\leq B\leq B_2$,
 by going to subsequence if necessary.
Because $|A_{\epsilon}|^{-\frac{1}{2}}:E\to E$ is compact,
$|A_{\epsilon}|^{-\frac{1}{2}}B_{\epsilon}(u_{n})|A_{\epsilon}|^{-\frac{1}{2}}
y_{n}\to
|A_{\epsilon}|^{-\frac{1}{2}}B_{\epsilon}|A_{\epsilon}|^{-\frac{1}{2}}
y$, and from (5.2) it follows that $y_n\to y$ and
  $$
y^{+}-y^{0}-y^-
-|A_{\epsilon}|^{-\frac{1}{2}}B_{\epsilon}|A_{\epsilon}|^{-\frac{1}{2}}
y=0
 $$
Set $x=|A_{\epsilon}|^{-\frac{1}{2}}y$. Then $x\neq 0$ since
$||y||=1$, and $ Ax-Bx=0$. This is impossible because Proposition
1.5(iv) implies that $\nu _{A}(B)=0$. And the proof is
complete.\hfill\hb

\vskip4mm

In order to prove Proposition 5.2 we need the following two lemmas.

 {\bf Lemma 5.3}. Suppose assumptions (i-ii) in Theorem 1.1 hold. Then there exists $R_{0}>0$ such that
\begin{eqnarray}
<a'(u_{1}+u_{2}),u_{2}-u_{1} >1 \hbox{ as } \|u_{2}\|\geq R_{0},
\hbox{ or } \|u_{1}\|\geq R_{0},
\end{eqnarray}
where $u_{2}\in E^{+}_{A,\beta}(B_{2})$ and $u_{1}\in
E^{-}_{A,\beta}(B_1)$.

{\bf Proof}. For any $u^{+}+u^0=u_{1}+u_{2}$ with $u_{2}\in
E^{+}_{A,\beta}(B_{2})$ and $u_{1}\in E^{-}_{A,\beta}(B_{1})$,
from (3.12) and assumption (ii) we have
  \begin{eqnarray}
<a'(u^{+}+u^{0}),u_{2}-u_{1}
>=&&(u^{+}-u^{0}
   -(P^{+}+P^{0})|A_{\epsilon}|^{-\frac{1}{2}}(B_{\epsilon}(|A_{\epsilon}|^{-\frac{1}{2}}u)
                                                   |A_{\epsilon}|^{-\frac{1}{2}}u+C(|A_{\epsilon}|^{-\frac{1}{2}}u)),
u_{2}-u_{1})\nn\\
=&&((P^{+}-P^{0})u_{2}-(P^{+}+P^{0})B_{\epsilon}(|A_{\epsilon}|^{-\frac{1}{2}}u)|A_{\epsilon}|^{-\frac{1}{2}}u_{2},
u_{2})\nn\\
&&-((P^{+}-P^{0})u_{1}-(P^{+}+P^{0})B_{\epsilon}(|A_{\epsilon}|^{-\frac{1}{2}}u)|A_{\epsilon}|^{-\frac{1}{2}}u_{1},
u_{1})+r(u_1,u_2,\beta)\nn\\
\geq&&
\bar{q}_{A,\beta;B_2}(u_{2},u_{2})-\bar{q}_{A,\beta;B_1}(u_{1},u_{1})+r(u_1,u_2,\beta),
   \end{eqnarray}
where $u=u_2+u_1+u^-$ and
  \begin{eqnarray}
u^{-}&&=-P^{-}|A_{\epsilon}|^{-\frac{1}{2}}\Phi_{\epsilon}'(u)(|A_{\epsilon}|^{-\frac{1}{2}}u)
         =-P^{-}|A_{\epsilon}|^{-\frac{1}{2}}(B_{\epsilon}(|A_{\epsilon}|^{-\frac{1}{2}}u)|A_{\epsilon}|^{-\frac{1}{2}}u
                                                       +C(|A_{\epsilon}|^{-\frac{1}{2}}u))\nn\\
    r(u_1,u_2,\beta)&&=-((P^{+}+P_0)|A_{\epsilon}|^{-\frac{1}{2}}(B_{\epsilon}(|A_{\epsilon}|^{-\frac{1}{2}}u)
                                        |A_{\epsilon}|^{-\frac{1}{2}}u^-+C(|A_{\epsilon}|^{-\frac{1}{2}}u)),u_2-u_1).\nn
  \end{eqnarray}
 From assumption (ii), let $M>0$ such that $\|B_{\epsilon}(x) \|\leq
 M$, $\|C(x) \|\leq M$ $\forall x\in X$. A simple calculation shows that
    \begin{eqnarray}
||u^{-}||\leq
\frac{\sqrt{\beta}M}{(\beta-M)\sqrt{\epsilon}}\|u^{\ast}\|+\frac{M^2
\sqrt{\beta}}{\beta-M},
  \end{eqnarray}
and
  \begin{eqnarray}
|r(u_1,u_2,\beta)|\leq M (\frac{M
}{(\beta-M)\epsilon}(\|u_{1}\|+\|u_{2}\|)+1+\frac{M
\sqrt{\beta}}{\beta-M})(\|u_{1}\|+\|u_{2}\|).\nn
  \end{eqnarray}
Then (5.3) follows from Proposition 4.3(i) and (5.4). The proof is
complete.\hfill\hb

{\bf Lemma 5.4}. Let $R_{0}$ be defined in Lemma 5.3. Then
$a(u_{2}+u_{1})\rightarrow -\infty$ uniformly for $u_{2}\in
E^{+}_{A,\beta}(B_{2})\cap B_{R_{0}}$ as $\|u_{1}\|\rightarrow
+\infty$.

{\bf Proof}. Because
  \begin{eqnarray}
&&\Phi(|A_{\epsilon}|^{-\frac{1}{2}}u)
           =\int_0^1(\Phi_{\epsilon}'(|A_{\epsilon}|^{-\frac{1}{2}}u\theta),|A_{\epsilon}|^{-\frac{1}{2}}u)d\theta\nn\\
  &&=\int_0^1(B_{\epsilon}(|A_{\epsilon}|^{-\frac{1}{2}}u\theta)|A_{\epsilon}|^{-\frac{1}{2}}u\theta
                                     +C(|A_{\epsilon}|^{-\frac{1}{2}}u\theta),|A_{\epsilon}|^{-\frac{1}{2}}u)d\theta\nn\\
  &&\geq{1\over
  2}((B_{1})_{\epsilon}|A_{\epsilon}|^{-\frac{1}{2}}u,|A_{\epsilon}|^{-\frac{1}{2}}u)-{M\over\sqrt{\epsilon}}||u||.
   \end{eqnarray}
By assumption (ii), for any $u^0+u^+=u_1+u_2$ with $u_2\in
E^+_{A,\beta}(B_{2})\cap B_{R_0}, u_1\in E^-_{A,\beta}(B_{1})$,
   \begin{eqnarray}
 a(u_{2}+u_{1})&&=\frac{1}{2}\|u^{+}\|^{2}-\frac{1}{2}\|u^{0}\|^{2}-
    \frac{1}{2}\|u^-\|^2-\Phi_{\epsilon}({|A_{\epsilon}|^{-\frac{1}{2}}u})\nn\\
               &&\leq \frac{1}{2}{\bar
               q}_{A,\beta;B_1}(u_1,u_1)+\frac{C_{2}}{\sqrt{\beta}}\|
               u_{1}\|^2+C_{3}\nn
                \end{eqnarray}
for $\beta$ larger enough and some constants $C_{2}>0$, $C_{3}>0$.
Here we used (5.4) and (5.5). By (4.11) (4.13) and (4.14)
$\frac{1}{4}{\bar
q}_{A,\beta;B_1}(u_1,u_1)+\frac{C_{2}}{\sqrt{\beta}}\|
               u_{1}\|^2 \leq 0$ for $\beta$ larger enough.
 Thus, $a(u_{2}+u_{1})\rightarrow
-\infty$ as $\|u_{1} \|\rightarrow \infty$ uniformly for $u_2\in
E^+_{A,\beta}(B_{2})$ for $\beta>0$ larger enough. The proof is
complete.\hfill\hb

\vskip4mm We will use Lemma 3.2 again to investigate critical
points of the functional $a(u^++u^0)$ and we need to calculate
some relative homology groups. As in Chang[5] and in
Mawhin-Willem[16] we say that the topological space pair $(X',Y')$
with $X'\subset Y'$  is the deformation retract of a topological
space pair $(X,Y)$ with $Y\subset X$ if $X'\subset X$, $Y'\subset
Y$ and there exists $\eta:[0,1]\times X\to X$ satisfying
     $$
     \eta(0,\cdot)=id_X, \eta(1,X)\subset X',\eta(1,Y)\subset Y,
     \eta(t,Y)\subset Y,
     $$
and
   $$
\eta(t,\cdot)=id_{X'},   \forall t\in[0,1].
    $$
It is well-known that if $(X',Y')$ is a deformation retract of
$(X,Y)$, then
   $$
   H_q(X,Y;\R)\cong H_q(X',Y';\R), q=0,1,2,\cdots.
   $$
For any $X\subset Y\subset Z$ if there exists $\tau: [0,1]\times
Y\to Y$ satisfying $\tau(0,\cdot)=id_Y, \tau(1,Y)\subset Z$ and
$\tau(t,\cdot)_Z=id_Z$, then $Z$ is called a strong deformation
retract of $Y$. And from a result in [16, page 171] by
Mawhin-Willem we have
   $$
H_q(X,Y;\R)\cong H_q(X,Z;\R), q=0,1,2,\cdots.
   $$

{\bf Proof of Proposition 5.2} Set ${\cal
M}_{R_0}=(E^+_{A,\beta}(B_2)\bigcap B_{R_0})\bigoplus
E^-_{A,\beta}(B_1), \sigma(t,u)=e^{-t}u_2+e^tu_1$ and $T_u=\ln
||u_2||-\ln R_0$ if $u=u_2+u_1$ and $\|u_2\|>R_0$. Here $R_0$ is
defined in Lemma 5.3 and $B_R:=\{u\in E|\|u\|\leq R\}$. Define
 \begin{eqnarray}
 \eta(t,u_2+u_1)&&=u_2+u_1,||u_2||\leq R_0,\nn\\
                &&=\sigma(T_ut,u), ||u_2||>R_0.\nn
 \end{eqnarray}
By Lemma 5.3 it is easy to verify that $({\cal M}_{R_0},{\cal
M}_{R_0}\bigcap a_c)$ is a deformation retract of $(E, a_c)$ for any
$c\in\R$. And hence,
  \be
  H_q(E,a_c;\R)=H_q({\cal M}_{R_0},{\cal M}_{R_0}\bigcap a_c; \R),
  q=0,1,2,\cdots.\nn
  \ee
Now we begin to prove  that
    \begin{eqnarray}
   H_q({\cal M}_{R_0},{\cal M}_{R_0}\bigcap
   a_{-c};\R)\cong\delta_{q\gamma}\R,\nn
        \end{eqnarray}
where $-c<a(\theta)$. By Lemma 5.4, there exist $T>0, c_1>c_2>T,
R_1>R_2>R_0$ such that
     \begin{eqnarray}
{\cal N}_{R_1}\subset a_{-c_1}\bigcap{\cal M}_{R_0}\subset{\cal
N}_{R_2}\subset a_{-c_2}\bigcap{\cal M}_{R_0}\subset{\cal
N}_{R_0},\nn
   \end{eqnarray}
where ${\cal N}_R:=(E^+_{A,\beta}(B_2)\bigcap
B_{R_0})\bigoplus(E^-_{A,\beta}(B_1)\backslash B_R)$ for any
$R>0$. For any
 $u\in {\cal M}_{R_0}\bigcap(a_{-c_2}\setminus a_{-c_1})$, since
$\sigma(t,u)=e^{-t}u_2+e^tu_1, a(\sigma(t,u))$ is continuous with
respect to $t$ and $a(\sigma(0,u))=a(u)\in(-c_1,-c_2]$. From (2.21)
    \begin{eqnarray}
 &&{d\over dt}a(\sigma(t,u))=\langle a'(\sigma(t,u)),\sigma'(t,u)\rangle\nn\\
 &&=\langle a'(e^{-t}u_2+e^tu_1),-e^{-t}u_2+e^tu_1\rangle\leq-1\nn
      \end{eqnarray}
as $t>0$.  So the time $t=T_1(u)$ arriving at $a_{-c_1}\bigcap{\cal
M}_{R_0}$ exists uniquely and is defined by $a(\sigma(t,u))=-c_1$.
The continuity of $t=T_1(u)$ comes from the implicit function
theorem. Define
      \begin{eqnarray}
 \eta _1(t,u) &&=u, \ \ \ x\in a_{-c_1}\bigcap {\cal M}\nn\\
              &&=\sigma(T_1(u)t,u), \ \ \ u\in {\cal M}\bigcap (a_{-c_2}\backslash a_{-c_1});\nn
     \end{eqnarray}
and
    \begin{eqnarray}
 \eta_2 (t,u)&&=u,  \|u_1\| \geq R_1\nn\\
           &&=u_2+tu_1+(1-t)\frac{u_1}{\|u_1\|}R_1, \ \|x_1\|<R_1.\nn
     \end{eqnarray}
By the map $\eta(t,u)=\eta_2(t,\eta_1(t,u))$ we can verify that
$(E^+_{A,\beta}(B_2)\bigcap B_{R_0}) \oplus
(E^-_{A,\beta}(B_1)\backslash int (B_{R_1})$ is  a strong
deformation retract of ${\cal M}\bigcap a_{-c_2}$:
  \begin{eqnarray}
&&\eta(t,u)=u\ \ \forall\ \ u\in (E^+_{A,\beta}(B_2)\bigcap B_{R_0})
  \oplus (E^-_{A,\beta}(B_1)\backslash int (B_{R_1}), t\in
      [0,1],\nn\\
&&\eta(0,u)=u\ \ and\ \ \eta(1,u)\in (E^+_{A,\beta}(B_2)\bigcap
B_{R_0}) \oplus (E^-_{A,\beta}(B_1)\backslash int (B_{R_1}), \ \
\forall\ \ u\in {\cal M}\bigcap a_{-c_2}.\nn
   \end{eqnarray}
From a result in [16, page 171] it follows that
   \begin{eqnarray}
    &&H_q({\cal M},{\cal M}\bigcap a_{-c_2};\R)\nn\\
        &&\ \ \cong  H_q((E^+_{A,\beta}(B_2)\bigcap B_{R_0})
             \bigoplus E^-_{A,\beta}(B_1),(E^+_{A,\beta}(B_2)\bigcap B_{R_0}) \bigoplus
(E^-_{A,\beta}(B_1)\backslash int (B_{R_1})) ;\R)\nn\\
        &&\ \ \cong H_q(E^-_{A,\beta}(B_1)\bigcap B_{R_1},\partial (E^-_{A,\beta}(B_1)\bigcap
                              B_{R_1});\R)\nn\\
        &&\ \ \cong \delta_{q\gamma}\R.\nn
   \end{eqnarray}
Therefore, combining (5.7) implies that  (5.1) holds and the proof
is complete.\hfill\hb

 {\bf Proof of Theorem 1.1.} From Proposition 4.2 (iii), Definition 1.3
 and Theorem 1.5(ii)
\begin{eqnarray}
i_{A,\beta}(B_{1})-m^{-}(a''(\theta))=i_{A,\beta}(B_{1})-i_{A,\beta}(B_{0})=i_{A}(B_{1})-i_{A}(B_{0}).\nn
\end{eqnarray}
So assumption (iii) implies that $\gamma=i_{A,\beta}(B_{1})\notin
[i_{A,\beta}(B_{0}), i_{A,\beta}(B_{0})+\nu_{A,\beta}(B_{0})]$. And
$\nu_{A,\beta}(B_{0})=0$ means that $\theta$ is a non-degenerate
critical point; $\nu_{A}(\Phi''(x_{0}))\leq
|i_{A}(B_{1})-i_{A}(B_{0}) |$ implies that $m^{0}(a''(x_{0}))\leq
|i_{A,\beta}(B_{1})-i_{A,\beta}(B_{0})|$. By Lemma 3.2 and
Propositions 5.1 and 5.2, the proof is complete.\hfill\hb

{\bibliographystyle{abbrv}

\end{document}